\documentclass[preprint,12pt,3p]{elsarticle}

\usepackage{amssymb}
\usepackage{bm}
\usepackage{amsmath}
\usepackage{float}
\usepackage{comment}
\usepackage{algorithm}
\usepackage{algorithmic}
\usepackage{diagbox}
 \usepackage{graphicx}
  \usepackage{caption}
 \usepackage{booktabs}
 \usepackage{color}
 \usepackage{verbatim}
 \usepackage{subcaption}
 \usepackage{ulem}
\usepackage{overpic}
\usepackage{rotating}

\def\S{\mathcal{S}}

\usepackage{lineno}

\journal{\ \ }

\begin{document}


\begin{frontmatter}

\title{Learning PDEs from data on closed surfaces with sparse optimization}
\tnotetext[label0]{Funding: The work of the first author was supported  by  NSFC (No. 12101310), NSF of Jiangsu Province (No. BK20210315), and the Fundamental Research Funds for the Central Universities (No. 30923010912), the work of the second author was supported by the General Research Fund (GRF No. 12301520, 12301021, 12300922) of Hong Kong Research Grant Council, and the work of the corresponding  author was supported by NSFC(No. 11901377) and NSFC(No. 12171093).}

\author[label1]{Zhengjie Sun}
\ead{zhengjiesun@njust.edu.cn}

\author[label3]{Leevan Ling}
\ead{lling@hkbu.edu.hk}

\author[label2]{Ran Zhang\corref{cor1}}

\cortext[cor1]{Corresponding author}
\address[label2]{School of Mathematics, Shanghai University of Finance and Economics, Shanghai, China.}
\ead{zhang.ran@mail.shufe.edu.cn}

\address[label1]{School of Mathematics and Statistics, Nanjing University of Science and Technology, Nanjing, China.}

\address[label3]{Department of Mathematics, Hong Kong Baptist University, Kowloon Tong, Hong Kong.}





\begin{abstract}
Discovering underlying partial differential equations (PDEs) from observational data has important implications across fields. It bridges the gap between theory and observation, enhancing our understanding of complex systems in applications.
In this paper, we propose a novel approach, termed physics-informed sparse optimization (PIS), for learning surface PDEs. Our approach incorporates both $L_2$ physics-informed model loss and $L_1$ regularization penalty terms in the loss function, enabling the identification of specific physical terms within the surface PDEs. 
The unknown function and the differential operators on surfaces are approximated 
by some extrinsic meshless methods. 
We provide practical demonstrations of the algorithms including linear and nonlinear systems. 
The numerical experiments on spheres and various other surfaces demonstrate the effectiveness of the proposed approach in simultaneously achieving precise solution prediction and identification of unknown PDEs. 

\end{abstract}


\end{frontmatter}



\section{Introduction}

Data-driven modeling only with available data has been widely considered in learning theory and its associated application areas. While many phenomena in science and engineering can be formulated as partial differential equations (PDEs), traditional PDE models predominantly rely on system behavior descriptions and classical physical laws. 
Consequently, the fundamental challenge lies in distilling the underlying PDEs from the given data. 

The methodology of data-driven modeling can even date back to the time of Kepler, who employed the most meticulously guarded astronomical data of the day to discover a data-driven model for planetary motion. Later in 1975, Gauss introduced the least squares regression (LSR) algorithm, which provided a numerical framework for learning underlying models from data. 
The classical Prony's method was originally developed to use the difference equation as a discrete analogue of the linear ordinary differential equation, then used LSR to compute coefficients of linear ordinary difference equations. However, it is known to perform poorly in the presence of noisy samples. To address this limitation, various stabilization and modification methods for Prony's method have been proposed, as discussed in Osborne et al. \cite{Osborne95} and Zhang et al. \cite{Zhang}, for example. 
Furthermore, Schmidt and Lipson in \cite{Science} introduced symbolic regression and evolutionary algorithms to directly learn physical laws, such as nonlinear energy conservation laws and Newtonian force laws, from experimentally captured data.  

With the rapid advancement of data storage and data science tools, data-driven discovery of potential models has entered a new era with the emergence of big data. 
Numerical simulations including Lorentz system (ODE) and fluid flow (PDE) are presented in \cite{Brunton} with thresholded least square method to promote sparsity.
A deep PDE network has been considered to deal with the multi-dimensional systems based on the convolution kernels in \cite{Dong}.
\cite{PINN,Shekarpaz} proposed the physics-informed neural networks (PINN) framework to address both data-driven solutions and data-driven discovery of PDEs. 
Wu et al. \cite{Wu19} sucessfully applied the sparse regression to learn a chaotic system, effectively simulating the sustaining oscillations during the sedimentation of a sphere through non-Newtonian fluid.  Wu et al. \cite{Wu_AA} further discussed two probabilistic solutions for the dynamical system, namely the random branch selection iteration (RBSI) and random switching iteration (RSI). Other related works include, for instance, parameter identifications \cite{GaoZhen} in uncertainty quantification  and calibration \cite{Hildebrandt} in the financial industry, which often involve more prior knowledge about the corresponding PDE models. 

In recent years, surface PDEs find a wide range of applications in various fields, including imaging processes, biological reactions, fluid dynamics, and computer graphics \cite{Adalsteinsson, Auer, Biddle, Charles}. As a result, there has been significant research on numerical methods to solve surface PDEs, employing various techniques such as intrinsic, extrinsic, and embedding methods. Intrinsic methods involve local parameterization of surfaces and discretization of surface differential operators on surface meshes \cite{Floater}. Extrinsic methods transform the differential operators on surfaces into extrinsic coordinates \cite{chen}, while embedding approaches extend the surface PDEs to embedding spaces \cite{Ruuth}. Although there have been advancements, such as the use of physics-informed convolutional neural networks (PICNN) to solve PDEs on spheres \cite{Lei}, the inverse problem of discovering hidden PDEs based on surface data remains an open challenge. In this paper, our focus is on addressing this challenge by developing novel techniques for the discovery of nonlinear PDEs on closed surfaces. 

Let  $\mathcal{S} \subset \mathbb{R}^d$ be a closed smooth surface,  $d_{\mathcal{S}}$ denotes the dimension, and $d_{co}:=d-d_\mathcal{S}$ is the codimension. Suppose that  $u=u(\bm{x},t)$ 
satisfies the following evolutionary PDE on the surface $\mathcal{S}$,

\begin{equation}\label{t_sPDE}
\left\{
\begin{aligned}
    &\frac{\partial u}{\partial t}=\mathcal{L}_{\mathcal{S}}u-f, \ \ \bm{x} \in \mathcal{S},\ \  t \in (0,T], \\
    &u(\bm{x},0)=u_0(\bm{x}), 
\end{aligned}
\right.
\end{equation}
where $\mathcal{L}_{\mathcal{S}}$ is the unknown differential operator defined on $\mathcal{S}$.

The differential operators on the surface can be defined in terms of the standard Euclidean differential operators. For any $\bm{p} \in \mathcal{S}$, let $\bm{n}=\bm{n}(\bm{p})=(n_1(\bm{p}),n_2(\bm{p}),\dots,n_d(\bm{p}))^\top$ be the unit outward normal vector at $\bm{p}$ and $\mathcal{P}$ be the corresponding orthogonal projection operator from $\mathbb{R}^d$  onto the tangent space $\mathcal{T}_{\bm{p}}\mathcal{S}$. Then we can define the surface gradient and surface Laplace-Beltrami operator as
\begin{equation}\label{grad_S} 
\nabla_{\mathcal{S}}:= \mathcal{P}\nabla   = (I_d-\bm{n}\bm{n}^\top)\nabla \ \ \ \text{and}\ \ \ \Delta_{\mathcal{S}}:=  \nabla_{\mathcal{S}} \cdotp \nabla_{\mathcal{S}}
\end{equation}
with $I_d$ being the identity operator of size $d \times d$. The surface gradient and Laplacian-Beltrami operator can be analytically transformed to  Euclidean forms by using the equation \eqref{grad_S}. Thus, we could approximate the underlying function $u$ and the differential operators in $\mathcal{L}_{\mathcal{S}}u$ by using the extrinsic 
method with the meshfree radial basis function (RBF). Moreover, as a fully data-driven approach, if we are only provided with a point cloud and do not have the analytical expression for the surface, the normal extension method \cite{OGM} can still be employed to approximate the surfaces' normal vectors.  

We exploit the understanding that the majority of physical systems are governed by a limited number of nonlinear terms. 
In the process of learning the nonlinear unknown $\mathcal{L}_{\mathcal{S}}u$, the sparse optimization will be proposed to minimize the model loss function. 
Except the $L_2$ physics-informed model loss, we also introduce the $L_1$ regularization penalty term to identify the physical terms in the PDE. 
This enables us to achieve a balance between model complexity and data fidelity, leading to a more parsimonious and accurate representation of the underlying phenomenon. 

The novelty of this paper can be attributed to two key aspects. Firstly, in addition to identifying the unknown parameters in the model, which is a common focus in existing works, the paper also determines terms with a clear physical background, such as diffusion or convection terms. Secondly, our method stands out by eliminating the need for iterations and offering solutions that surpass the accuracy of the PINN for surface PDEs, as demonstrated in \cite{Fu}.

The paper is structured as follows: In Section 2, we introduce the methodology for learning stationary surface PDEs. Subsequently, in Section 3, we generalize the approach to evolutionary surface PDEs. 
Finally, in Section 4, we present numerical simulations that demonstrate the capability of our method to identify underlying PDE models and predict long-term solutions on various closed surfaces.

\section{Discovery of stationary PDEs on surfaces}
\label{sec1}

In this section, we focus on the stationary equation of \eqref{t_sPDE} defined as
\begin{equation}\label{sPDE}
\mathcal{L}_{\mathcal{S}}u(\bm{x})=f, \ \ \bm{x} \in \mathcal{S}.
\end{equation}
To learn this stationary PDE (\ref{sPDE}), we are provided with a set of samples $\{\bm{x}_i, u^*(\bm{x}_i)\}_{i=1}^{N}$. Here 
the distinct nodes $X=\{\bm{x}_i\}_{i=1}^{N} \subset  \mathcal{S} \subset   \mathbb{R}^d$ are exactly located on the surface $\mathcal{S}$, and the function value 
 $u^*(X)= u^*|_X
 $ may contain noise. 
Note that the presence of sampling noises may cause the nodes in $X$ to deviate slightly from the actual surface $\mathcal{S}$. In such situations, approaches like the closest point method, discussed in \cite{cheung}, can be employed. However, for the purposes of this paper, we do not consider such cases.


We propose a three-step process to discover the hidden stationary PDEs on surfaces. In the first step, we construct a model library of potential functions and surface differential operators using a polynomial-based approach. This library serves as a foundation for capturing the underlying dynamics of the system. In the second step, we provide numerical approximations of the function and surface differential operators using the meshless kernel-based method. Finally, in the third step, we employ sparse optimization techniques to uncover the underlying model. The goal is to find a sparse representation of the constructed model, highlighting the essential terms and features that contribute significantly to the system's behavior.

\subsection{Model library construction}

In reality, the operator $\mathcal{L}_{\mathcal{S}}$ may include higher-order nonlinear differential operators in the hidden PDE. However, when dealing with limited observations of the quantity $u$ on the surface nodes, we need to take into account both the physical background and computational complexity. To strike a balance between accuracy and feasibility, we consider the operator $\mathcal{L}_{\mathcal{S}}u := \mathcal{L}_{\mathcal{S}}(u,\nabla_{\mathcal{S}}u,\Delta_{\mathcal{S}}u)$ up to the second order because second-order PDEs have been extensively studied and often provide sufficient accuracy to describe various physical phenomena.

To approximate the unknown operator $\mathcal{L}_{\mathcal{S}}u$, we use a polynomial-based approach to  establish a library of candidate nonlinear differential operators on the surface up to second order. First, let $\{\bm{e}_k(\bm{x})\}_{k=1}^d$ denotes the orthogonal moving frame at $\bm{x} \in \mathcal{S}$.  We define a map $\bm{\lambda}: \mathcal{S} \rightarrow \mathbb{R}^{d+2}$  as

\begin{equation}\label{z}
\bm{z}:=\bm{\lambda}(\bm{x})=(u(\bm{x}),
[\nabla_{\mathcal{S}} u]_1(\bm{x}),\ldots,
[\nabla_{\mathcal{S}} u]_d(\bm{x}),\Delta_{\S}u(\bm{x})),
\end{equation}
where $u$ is the unknown function, 
$[\nabla_{\mathcal{S}} u]_k(\bm{x}):=\bm{e}_k \cdot \nabla_{\mathcal{S}} u(\bm{x})$ is the $k$th ($1\leq k \leq d$) component of $\nabla_{\mathcal{S}} u$, 
and $\Delta_{\mathcal{S}}u$ is the surface Laplacian. With the map $\bm{\lambda}$, we introduce multivariate monomials in $\mathbb{R}^{d+2}$ given by

$$
p_{\bm{\alpha}}(z)=\bm{z}^{\bm{\alpha}}=z_1^{\alpha_1}z_2^{\alpha_2}\cdots z_{d+2}^{\alpha_{d+2}},~~\bm{\alpha}=(\alpha_1,\ldots,\alpha_{d+2}).
$$

Next, we construct the model library, denoted by $\Lambda$, consisting of all polynomials with degrees less than $\ell$. It can be expressed as 
\begin{equation}  \Lambda:=\Big[p_{\bm{\alpha}}\Big],~~|\bm{\alpha}|\leq \ell,
\end{equation}
where $\ell$ is the highest degree of polynomials included in the library.

Let $\bm{\xi} = [\xi_1, \xi_2, \dots,\xi_n]^\top$ with $n=\binom{d+2+\ell}{\ell}$ be the coefficient vector of the polynomial bases in the library $\Lambda$, then we can construct the approximation $\Lambda\bm{\xi}$ to $\mathcal{L}_{\mathcal{S}}(u,\nabla_{\mathcal{S}}u,\Delta_{\mathcal{S}}u)$ and solve the following equation 
\begin{equation}
\Lambda\bm{\xi}=f
\end{equation}
This equation holds for all the samples, allowing us to define the feature matrix 
\begin{equation}\label{library}
\Lambda(X)=
\left[p_{\bm{\alpha}}|_X\right],
\end{equation} 
where $p_{\bm{\alpha}}|_X$ represents the polynomial $p_{\bm{\alpha}}$ features evaluated at the sample points $X$.

Finally, we can approximate the stationary PDE by solving the following linear system  
\begin{equation}
\Lambda(X)\bm{\xi}= f(X),
\end{equation}
with $f(X)=f|_X$. It's worth noting that although the bases in the library $\Lambda$ contain many nonlinear operators, we only need to deal with a linear problem regarding the coefficients of these bases.

Dynamical systems that capture all relevant spatial scales can have enormous dimensions, as stated by Temam \cite{Temam}. It is unlikely for any specific model to perfectly simulate the true dynamics, and improvement can only be achieved progressively. Therefore, it is possible to consider higher-order nonlinear differential equations or include higher-degree polynomials and trigonometric terms, such as $u^4$ and $\sin{u}$, to account for more complex dynamics. However, it is essential to strike a balance between incorporating additional terms and the associated physical backgrounds and computational complexity. Second-order PDEs of position have been widely studied in physics to describe various phenomena.  The higher order differential equations will be investigated in the Example 4 of Section 4.

Moreover, choosing the appropriate basis functions to sparsely represent the dynamics of real-world systems can be a challenging task. 
In Example 3 of Section 4, a different strategy for constructing the feature library is presented, which does not explicitly include the true differential operators in the equation. 

In mathematics and physics, simpler models that can reveal and interpret the underlying laws of physics are preferred. In many systems, $\mathcal{L}_{\mathcal{S}}u$ only contains a few dominant terms, making it sparse in specific operator spaces. Consequently, we expect only a few features to be activated. To identify the true physical terms in the surface PDE, we propose utilizing $L_1$ sparse optimization, which penalizes the sum of the absolute values of the model coefficients.
Most existing works have focused on identifying the unknown parameters associated with specific terms in the model. However, we can go beyond that and determine terms with physical backgrounds, such as diffusion or convection terms, using our approach. By incorporating prior knowledge of the underlying physics, we can enhance the interpretability of the model.



\subsection{Data interpolation and estimation of the  differential operators on  surfaces}

In our second step of data-driven modeling, we provide the numerical approximation of the function $u$ and the surface differential operators in \eqref{library} from the samples. 

In this paper, we focus on the radial basis function (RBF), a kernel that depends only on the distance between its two arguments. The advantage of RBF kernels is that they are meshless, meaning they do not require equidistant samples like the methods based on Prony's related methods.
RBF kernels have been successfully applied to develop the approximation of the unknown function and its differential operators on various manifolds \cite{Fuselier,Gia}. 

We start from the global kernels $\Phi_{\tau}: \mathbb{R}^d \times \mathbb{R}^d \rightarrow \mathbb{R}$ that are symmetric positive definite with Fourier transforms $\widehat{\Phi}_{\tau}$ satisfying the decay
$$
c_1(1+||\omega||_2^2)^{-\tau} \leq \widehat{\Phi}_{\tau}(\omega) \leq c_2(1+||\omega||_2^2)^{-\tau} \ \ \  \ \ \mathrm{for \ all} \ \omega \in \mathbb{R}^d
$$
for some constants $0<c_1 \leq c_2$.  One typical example is the  standard Whittle-Mat\'{e}rn-Sobolev kernels 
\begin{equation}\label{kernel_sob}
\Phi_{\tau}(x,y):=||x-y||_2^{\tau-d/2}\mathcal{K}_{\tau-d/2}(||x-y||_2) 
\end{equation}
where $\mathcal{K}$ is the modified Bessel function of the second kind. Another example is the class of Wendland compactly supported kernels. These $\Phi_{\tau}$ can reproduce  $\mathcal{H}^{\tau}(\mathbb{R}^d)$ 
for any integer $\tau > d/2$. We restrict this global kernels on $\mathcal{S}$ to obtain restricted kernels $\Psi_m: \mathcal{S} \times \mathcal{S} \rightarrow \mathbb{R}$, i.e., 
\begin{equation}\label{kernel_surface}
\Psi_m (\cdot,\cdot):= \Phi_{\tau}(\cdot,\cdot)_{|\mathcal{S} \times \mathcal{S}}
\end{equation}
which reproduce $\mathcal{H}^m(\mathcal{S})$  for  $m=\tau -d_{co}/2> d_{\mathcal{S}}/2$, see Theorem 5 in  \cite{Fuselier}. 

Then the interpolant of  $u$  with  the samples is given by
\begin{equation}\label{interpolant}
\widetilde{u}(\cdot)=\sum_{i=1}^N \kappa_i \Psi_m(\cdot,\bm{x}_i),
\end{equation}
where the coefficients $\bm{\kappa}=[\kappa_1, \kappa_2,\dots, \kappa_N]^\top$ can be obtained by solving the linear system $\Psi_m(X,X)\bm{\kappa}=u^*(X)$ with the entries given by $(\Psi_m(X,X))_{ij}=\Psi_m(\bm{x}_i,\bm{x}_j), i,j=1,\dots,N$. 
This Gram matrix is positive definite (so invertible) since the kernel is positive definite.
The analytical approximation error of this interpolant is guaranteed by \cite[Thm.~10]{Fuselier}. 
 
The surface differential operators need to be computed numerically. 
Under the framework of the extrinsic method, we have to transform the differential operators on surfaces to extrinsic coordinates via the projection $\mathcal{P}$ as shown in \eqref{grad_S}. 
In the scenario where only point clouds $X$ defining the surface $\mathcal{S}$ are available, identifying the implicit surface and obtaining the normal vector $\bm{n}$ of the surface become challenging.  In such cases, we use a RBF-based surface reconstruction technique developed in \cite{Turk} to approximate the underlying surface represented by the point cloud.   We employ the normal extension method proposed in \cite{OGM} using radial basis functions to approximate the normal vector.


Specifically, at each point $\bm{x}_i$ ($i=1,\ldots,N$), we compute a rough approximation of the normal vector $\bm{n}(\bm{x}_i)$ to the surface and define two new points, $\bm{x}_i-\delta\bm{n}(\bm{x}_i)$ and $\bm{x}_i+\delta\bm{n}(\bm{x}_i)$, which are located at a distance $\delta$ in the normal direction to the surface on either side. Afterwards, we define the distance function 
\begin{equation}\label{distance}
s(\bm{x})=\sum_{i=1}^{N} \alpha_i \Phi_{\tau}(\bm{x},\bm{x}_i)+\beta_i \Phi_{\tau}(\bm{x}, \bm{x}_i-\delta\bm{n}(\bm{x}_i))+\zeta_i \Phi_{\tau}(\bm{x}, \bm{x}_i+\delta\bm{n}(\bm{x}_i))
\end{equation}
The $0$-level surface corresponds to the approximation of the surface $\mathcal{S}$, while the $-1$-level and $1$-level surfaces correspond to the points inside and outside the surface, respectively. 
 The coefficients $\alpha_i$, $\beta_i$ and $\zeta_i$ can be obtained by the interpolation conditions $s(\bm{x}_i)=0$, and $s(\bm{x}_i-\delta\bm{n}(\bm{x}_i))=-1$ and $s(\bm{x}_i+\delta\bm{n}(\bm{x}_i))=1$.
It is important to note that the parameter $\delta$ influences the quality of the interpolant.
The coefficients in \eqref{distance} can be approximated using all $3N$ points. The surface $\mathcal{S}$ is defined as the zero-isosurface of the function $s(\bm{x})$. The approximated unit normal direction will be given by
$$
\widetilde{\bm{n}}=(s_1,s_2,\dots,s_d)/\sqrt{s_1^2+s_2^2+\dots+s_d^2}.$$
The projection $\mathcal{P}$ in the aforementioned equations will be replaced by 
$
\widetilde{\mathcal{P}} = I_d-\widetilde{\bm{n}}\widetilde{\bm{n}}^\top
$.

Now, suppose that the surface normal vector and the projection $\mathcal{P}$ are given analytically or numerically by the above approximation method,  
we can use \eqref{interpolant} to define the estimation of the surface gradient operator as 
\begin{equation}\label{gradient}
\widetilde{\nabla}_{\mathcal{S}} u:= \nabla_{\mathcal{S}}  \widetilde{u}=\mathcal{P}\nabla \widetilde{u}=\left[\begin{array}{c}
\mathcal{P}_1 \cdot \nabla \\
 \vdots\\
 \mathcal{P}_d \cdot \nabla 
 \end{array} \right] \widetilde{u},
\end{equation}
where $\mathcal{P}_k$ is the $k$th column of $\mathcal{P}$, whose entries are $ \mathcal{P}_{ij}$, $i,j=1,2,\dots,d$. 
Here the gradient operator acts on the radial basis function (\ref{interpolant}), then we have  the $k$th ($1\leq k \leq d$) component of $\widetilde{\nabla}_{\mathcal{S}} u$,  
\begin{equation}\label{gradient1}
[\widetilde{\nabla}_{\mathcal{S}} u]_k = \left( \mathcal{P}^\top_k[\nabla \Psi_m(\cdot,X)]\right)\left( [\Psi_m(X,X)]^{-1} u(X)\right) := G_k(\cdot,X)\bm{\kappa},
\end{equation}
where $G_k(\cdot,X): \mathcal{S} \rightarrow \mathbb{R}^{N\times N}$ is the matrix-valued function which can be evaluated at any given point on the surface.

A similar technique can be applied to approximate the Laplace-Beltrami operator. Following \cite{chen,Fuselier13}, we denote the interpolant of the surface gradient as $\widetilde{\widetilde{\nabla}_{\mathcal{S}} u } $, whose $k$th component is $\widetilde{[\widetilde{\nabla}_{\mathcal{S}} u ]}_k$. Then we obtain the approximation to the Laplace-Beltrami operator
\begin{equation}\label{laplace}
\widetilde{\Delta}_{\mathcal{S}} u := \widetilde{\nabla}_{\mathcal{S}} \cdot \widetilde{\nabla}_{\mathcal{S}} u= 
\nabla_{\mathcal{S}} \cdot \widetilde{\widetilde{\nabla}_{\mathcal{S}} u } 
= \mathcal{P}\nabla \cdot  
\left[\begin{array}{c}
\widetilde{[\widetilde{\nabla}_{\mathcal{S}} u ]}_1\\
\vdots\\
\widetilde{[\widetilde{\nabla}_{\mathcal{S}} u]}_d
\end{array} \right]. 
\end{equation}
The above operator 
(\ref{laplace}) can be written into the following matrix form,
\begin{equation}\label{laplace1}
\widetilde{\Delta}_{\mathcal{S}} u = (\mathcal{P}^\top_1\nabla, \dots,\mathcal{P}^\top_d\nabla) \cdot  
\left[\begin{array}{c}
[\Psi_m(\cdot,X)] [\Psi_m(X,X)]^{-1} [\widetilde{\nabla}_{\mathcal{S}} u]_1(X)\\
\vdots\\\relax
[\Psi_m(\cdot,X)][\Psi_m(X,X)]^{-1} [\widetilde{\nabla}_{\mathcal{S}} u]_d(X)
\end{array} \right].
\end{equation}
With the notations in \eqref{gradient1}, it reduces to 
\begin{equation}\label{laplace2}
\widetilde{\Delta}_{\mathcal{S}} u =\sum_{i=1}^N G_k(\cdot,X)[\Psi_m(X,X)]^{-1}G_k(X,X)\bm{\kappa}.
\end{equation}
Up to now, we have obtained all the discrete approximations to the function $u(x)$ and the surface differential operators by using radial basis function approximation methods, which could  be employed to construct the model library in \eqref{library}. Next, we shall use this library to discover the hidden PDE.

\subsection{Discovery of the differential equations}
Let's denote the approximation of the feature library as $\widetilde{\Lambda}(X)$, which is obtained by replacing $u(X)$, $\nabla_{\mathcal{S}} u(X)$ and $\Delta_{\mathcal{S}} u(X)$ in \eqref{library} with 
the estimations of the function value in \eqref{interpolant} and surface differential operators in \eqref{gradient} and \eqref{laplace2}, respectively. To determine the unknown coefficients, we solve the following linear problem
\begin{equation}\label{linear}
\widetilde{\mathcal{L}}_{\mathcal{S}}u(X) = f(X),
\end{equation}
where 
$
\widetilde{\mathcal{L}}_{\mathcal{S}}u(X) := \widetilde{\Lambda}(X)\bm{\xi}$.  

As mentioned earlier, models with more terms generally have smaller model loss, but there is a tradeoff between fidelity and generalization.
Therefore, in the following minimization problem, in addition to the squared $L_2$ norm of the model loss, we also introduce the $L_1$ regularization penalty term in LASSO \cite{Wu19} to determine a sparse vector of the coefficients
\begin{equation}\label{Lasso}
\widetilde{\bm{\xi}}=\mathop{\rm argmin}\limits_{\bm{\xi}} ||\widetilde{\mathcal{L}}_{\mathcal{S}}u(X) - f(X)||_2^2 + \mu||\bm{\xi}||_1,
\end{equation}
where $\mu > 0$ is a regularization parameter that controls the amount of shrinkage. The algorithm to compute \eqref{Lasso} involves solving a quadratic constrained optimization problem with linear constraints \cite{Schmidt}, 
\begin{equation}\label{constrained}
\mathop{\rm min}\limits_{\bm{\xi}} ||\widetilde{\mathcal{L}}_{\mathcal{S}}u(X) - f(X)||_2^2 + \mu\sum_{i=1}^{n}\gamma_i   \ \  s.t.   -\bm{\gamma} \leq \bm{\xi} \leq \bm{\gamma}.
\end{equation}
The $L_1$ regularization has been applied in various machine learning settings.  For example, in \cite{Tran}, the governing equations of a chaotic system are accurately recovered from corrupted data by solving a partial $L_1$ minimization problem with high probability. In \cite{Schaeffer}, the Douglas-Rachford algorithm is used to solve the $L_1$ least squares problem in learning partial differential equations. In the context of solving the multiscale elliptic PDE described in the reference \cite{Wang}, the approach employs $L_1$ regularization achieved by directly removing coefficients close to zero during specific iteration steps.

Besides, the square-root LASSO has also been considered recently as the sparsity-promoting optimization problem, which uses the $L_2$ norm of the model loss \cite{SRLasso}. 
It exhibits recovery performance analogous to that of LASSO, yet the optimal selection of its regularization parameter does not hinge upon the noise level present within the data. 
This work does not endeavor to comprehensively investigate the performance differential between the sparsity optimization methods. So we have merely included a cursory comparison of these two approaches in Example 1 of Section 4.

At the end of the section, we summarize the whole procedure in the following Algorithm 1. 

\begin{algorithm}
	\caption{Learning process of stationary PDEs}
	\begin{algorithmic} \label{Steady_PDE}
		\STATE \textbf{Initialization}: Given samples $\{\bm{x}_i, u^*(\bm{x}_i)\}_{i=1}^{N}$, and values of source function $\{f(\bm{x}_i)\}_{i=1}^{N}$ at the nodes on the surface $\mathcal{S}$.
		
		\STATE \textbf{Step 1}: Compute the interpolation to get the coefficients $\bm{\lambda}$ in (\ref{interpolant}).
		\STATE \textbf{Step 2}: 
		\IF{ normal vector $\bm{n}$ of  surface $\mathcal{S}$ is given}
		\STATE  
	    \ $\mathcal{P}$  is computed by theoretical normal vector of  $\mathcal{S}$; 
		\ELSE
		\STATE 
		 $\mathcal{S}$ is estimated as the zero-isosurface of $s(\bm{x})$ in (\ref{distance}), and the unit normal vector is computed by
	$\widetilde{\bm{n}}=(s_1,s_2,\dots,s_d)/\sqrt{s_1^2+s_2^2+\dots+s_d^2}$.
		\ENDIF
		\STATE \textbf{Step 3}: Compute $\widetilde{\nabla}_{\mathcal{S}} u$ in (\ref{gradient}) and $\widetilde{\Delta}_{\mathcal{S}} u $ in (\ref{laplace2}) to get the feature matrix $\widetilde{\Lambda}(X)$ in (\ref{linear}).
		\STATE \textbf{Step 4}: Solve (\ref{Lasso}) by $L_1$ regularization method (\ref{constrained}).
		\STATE  \textbf{Output}: $\widetilde{\bm{\xi}}$.
	\end{algorithmic}
	
\end{algorithm}


\section{Discovery of time-dependent PDEs on surfaces}

 For the underlying time-dependent PDEs on surface (\ref{t_sPDE}), 
the measurements now are provided as $\{(\bm{x}_i,t_j), u^*(\bm{x}_i,t_j)\}_{i=1, j=0}^{N,M}$ with $t_j=j\Delta t $. Denote $u^j(X)=u(X, t_j)$, 
$\nabla_{\mathcal{S}}u^j(X)=\nabla_{\mathcal{S}}u(X, t_j)$, and $\Delta_{\mathcal{S}}u^j(X)=\Delta_{\mathcal{S}}u(X,t_j)$, accordingly.
 
Similar to Eq. \eqref{z} to \eqref{library}, we define the map $\bm{\lambda}: \mathcal{S} \times \mathbb{R} \rightarrow \mathbb{R}^{d+2}$  as follows
\begin{equation}
\bm{z}:=\bm{\lambda}(\bm{x},t) =(u(\bm{x},t), [\nabla_{\mathcal{S}} u]_1(\bm{x},t)
\ldots, [\nabla_{\mathcal{S}} u]_d(\bm{x},t), 
\Delta_{\S}u(\bm{x},t)). 
\end{equation}
At the time level $t=t_j$, we construct the feature library of the differential operators as
\begin{equation}\label{library_t}
\Lambda^j(X)=
\Big[p_{\bm{\alpha}}^j|_X\Big], ~~|\bm{\alpha}|\leq \ell, 
\end{equation}
which contains polynomials with degrees less than $\ell$  of 
$(u, \nabla_{\mathcal{S}} u, \Delta_{\mathcal{S}} u)$ evaluated on the node set $X$ at time $t_j$. 
The equation \eqref{t_sPDE} holds for all time snapshots, which can be written as
\begin{equation}\label{eq:theequation}
\left.\frac{\partial u}{\partial t}(X)\right|_{t=t_j} = \Lambda^j(X)\bm{\xi}- f^j(X) \ \ j=0,1,\dots,M-1,
\end{equation}
where
$f^j(X)$ represents the values of $f$ evaluated on the nodes set $X$ at time $t_j$. The coefficients $\bm{\xi}$ are the unknowns to be determined.


Similar to the second step in Section \ref{sec1}, we  
approximate  $u^j(X)$ by using the radial basis function,
\begin{equation}\label{t_interpolant}
\widetilde{u}^j(\cdot)=\sum_{i=1}^N \kappa_i^j \Psi_m(\cdot,\bm{x}_i),
\end{equation}
where the coefficients $\bm{\kappa}^j=[\kappa_1^j, \kappa_2^j, \dots, \kappa_N^j]^\top$ are computed by solving the linear system $\Psi_m(X,X)\bm{\kappa}^j=u^{*,j}(X)$ at each time step.
Next we obtain the estimations of $\widetilde{\nabla}_{\mathcal{S}} u^j$ and $\widetilde{\Delta}_{\mathcal{S}} u^j$  from (\ref{gradient})-(\ref{laplace2}) accordingly. Furthermore, the approximation of  $\widetilde{\mathcal{L}}_{\mathcal{S}}u^j(X) := \widetilde{\Lambda}^j(X)\bm{\xi}$ can also be derived.  

For higher predictive accuracy, we can employ the second-order semi-implicit backward differentiation formula ($2-$SBDF) to discretize the equation \eqref{eq:theequation} in time. In this scheme, the diffusion term is discretized implicitly, resulting in the following equation 
{\small
\begin{equation} 
\label{2dt_sPDE}
\frac{ 3\widetilde{u}^{j+1}(X)-4\widetilde{u}^{j}(X)+\widetilde{u}^{j-1}(X)}{2\Delta t}=
\widetilde{\mathcal{L}}_{\mathcal{S}}^{\, 2}u^j(X)-(2f^j(X)-f^{j-1}(X)), \ \  j=1, \dots, M-1.
\end{equation} 
}
Here, $\Delta t$ represents the time step size. 
$\widetilde{\mathcal{L}}_{\mathcal{S}}^{\, 2}u^j(X)$ represents the second-order semi-implicit backward scheme, where $u(X)$, $\nabla_{\mathcal{S}}u(X)$, and $\Delta_{\mathcal{S}}u(X)$ are approximated by $2\widetilde{u}^j(X)-\widetilde{u}^{j-1}(X)$, 
$2\widetilde{\nabla}_{\mathcal{S}}u^{j}(X)-\widetilde{\nabla}_{\mathcal{S}}u^{j-1}(X)$, and 
$\widetilde{\Delta}_{\mathcal{S}}u^{j+1}(X)$, respectively, at each time step. 
Alternatively, the multi-time-stepping Runge-Kutta methods described in \cite{PINN} can also be employed for time discretization. 

Then we give the discretization approach with the aforementioned scheme (\ref{2dt_sPDE}). 
The mean squared model loss function with the regularization term is given by
{\small
\begin{equation}\label{t_Lasso}
\frac{1}{M}\sum_{j=0}^{M-1} \Big|\Big|
\frac{ 3\widetilde{u}^{j+1}(X)-4\widetilde{u}^{j}(X)+\widetilde{u}^{j-1}(X)}{2\Delta t}-
\widetilde{\mathcal{L}}_{\mathcal{S}}^{\, 2}u^j(X)+(2f^j(X)-f^{j-1}(X))
\Big|\Big|_2^2 +\mu||\bm{\xi}||_1.
\end{equation}
}
The corresponding algorithm is implemented in Algorithm 2.

\begin{algorithm}
	\caption{Learning process of  time-dependent PDEs}
	\begin{algorithmic} \label{}
			\STATE \textbf{Initialization}: Given samples $\{(\bm{x}_i,t_j), u^*(\bm{x}_i,t_j)\}_{i=1, j=0}^{N,M}$, time step  $\Delta t$ and value of  the right side function $\{f(\bm{x}_i,t_j)\}_{i=1, j=0}^{N,M}$ at the nodes on the surface $\mathcal{S}$.
			
			\STATE {Step 1}: 
			\IF{ the normal vector $\bm{n}$ of surface $\mathcal{S}$ is given}
			\STATE  
			\ $\mathcal{P}$ is computed by theoretical normal vector of  $\mathcal{S}$; 
			\ELSE
			\STATE 
			$\mathcal{S}$ is estimated as the zero-isosurface of $s(\bm{x})$ in (\ref{distance}), and the unit normal direction is
$\widetilde{\bm{n}}=(s_1,s_2,\dots,s_d)/\sqrt{s_1^2+s_2^2+\dots+s_d^2}$.
			\ENDIF
			\STATE {Step 2}: Compute  $G_k(X,X)$   and $G_k(X,X)[\Psi_m(X,X)]^{-1}G_k(X,X)$ in (\ref{gradient1}) and (\ref{laplace2}) respectively.
			\STATE {Step 3}: \FOR{$j=0$ to $M-1$}
			\STATE{Step 3.1:} Compute the interpolation to get the coefficients $\bm{\lambda}^j$ in (\ref{t_interpolant}); 
			 \STATE {Step 3.2}: Compute $\widetilde{\nabla}_{\mathcal{S}} u^j$ in (\ref{gradient}), and $\widetilde{\Delta}_{\mathcal{S}} u^j $ in (\ref{laplace2}), then get the feature matrix $\widetilde{\Lambda}^j(X)$ in (\ref{2dt_sPDE}).
			\ENDFOR
			
			\STATE {Step 4}: Solve the $L_1$ regularization (\ref{t_Lasso})  by quadratic constrained optimization.
			\STATE  \textbf{Output}: $\widetilde{\bm{\xi}}$.
	\end{algorithmic}
\end{algorithm} 

\section{Numerical demonstrations}

In this section, we will numerically test Algorithm 1 and Algorithm 2 by applying them to learn PDE models on various surfaces. The kernel functions involved in the following examples are standard Whittle-Mat\'{e}rn-Sobolev kernels in (\ref{kernel_sob}) restricted on the surfaces.  For time-dependent PDEs, the time discretization is performed using the second-order semi-implicit scheme $2-$SBDF, as described in (\ref{2dt_sPDE}).  In the case of an unit circle, the nodes are equally spaced, while for higher-dimensional surfaces, we employ algorithms from \cite{mesh} to generate the points. The noise imposed on the samples of the function values are independently and identically distributed Gaussian noise. 
In Example 1 and Example 2, the library of features defined in Eq.(\ref{library})  and 
Eq.(\ref{library_t}), respectively, consist of up to second-order polynomial terms involving $(u, \nabla_{\mathcal{S}} u, \Delta_{\mathcal{S}} u)$. 

In order to measure the accuracy, 
we introduce the relative $L_2$ error
$$
Re(t_j)=\frac{\sqrt{\sum_{i=1}^{N}(\widehat{u}(\bm{x}_i,t_j)-u(\bm{x}_i,t_j))^2}}{\sqrt{\sum_{i=1}^{N}(u(\bm{x}_i,t_j))^2}}
$$
where  $\widehat{u}(\bm{x}_i,t_j)$  is computed numerically by solving  the learned PDE using approximated Kansa method (AKM) \cite{chen} in Example 1 and Example 2, and $u(\bm{x}_i,t_j)$ is the true reference solution.

\subsection*{}\label{Ex1}

\noindent\textbf{Example 1: Linear stationary PDE}

\vspace{0.8em}

Consider the stationary PDE  
$$
\mathcal{L}_{\mathcal{S}}u=-\Delta_{\mathcal{S}}u+u=f, \ \ \ \  \bm{x} \in \mathcal{S}.
$$

\vspace{1em}



To begin with, we conduct a test on the unit circle using the samples $\{(x_i, y_i), u^*(x_i, y_i)\}_{i=1}^{N}$. In this case, we assume that the exact solution is given by $u=e^{x+y}(x^3+y^4+1)$, and we compute the corresponding function $f$ to satisfy the equation. 
The smoothness order $m$ of the kernels is chosen to be $6$ here.  

In the first experiment, where the samples are noiseless and the training dataset size is small ($N=30$), the errors for the coefficient term $\Delta_{\mathcal{S}}u$ and $u$ are found to be $0.007\%$ and $0.003\%$, respectively. However, the remaining terms in the library are not selected in the equation, indicating that they do not significantly contribute to the dynamics of the system. A regularization parameter of $\mu=0.01$ is used for this experiment with noiseless samples.

In the second experiment, a larger training dataset of size $N=300$ is used with a low noise level. In this case, the errors for the coefficient term $\Delta_{\mathcal{S}}u$ and $u$ are measured to be $0.51\%$ and $0.042\%$, respectively. Importantly, all the terms in the true equation are selected, suggesting that the model accurately captures the underlying dynamics of the system.
For this case of noisy samples, a higher regularization parameter of $\mu=20$ is chosen to mitigate the effects of noise and promote a smoother solution.

\begin{table} [H]
     \footnotesize
	\centering 
	\renewcommand\arraystretch{1.3}
	\begin{tabular}{ l|l 
		}
	\bottomrule[0.8pt]
		 \diagbox[width=10em]{Noise}{Model}  & True model  $-\Delta_{\mathcal{S}}u+u=f$ 
		\\
		\hline
		  $0\% (N=30)$& $ -1.0000  \Delta_{\mathcal{S}}u+1.0000 u=f $ 
		\\
		 $0.01\%(N=300)$&   $-0.9950  \Delta_{\mathcal{S}}u+0.9996 u=f $          
		\\ 

			\toprule[0.8pt] 
	\end{tabular}	
	\medskip
	
	\caption{
		{The identified model  in Example 1 on the unit circle with different number of samples and  noise levels.}}
	\label{Tab:ex1}
\end{table}

\vspace{1em}


\vspace{0.8em}

For the unit sphere, 
we test the exact solution $u=10xyz+5xy+z$. 
The smoothness order $m$ of the kernels is chosen as $4$. 
 To evaluate the accuracy of the normal vectors $\bm{n}$ computed using the analytical expression or the normal extension method, we examine the learned models presented in Table \ref{Tab:ex2_n} and Table \ref{Tab:ex2}, respectively. 

The results indicate that when noiseless samples are used with a smaller number of nodes, the normal extension method provides a slightly more accurate model compared to the analytical expression. For instance, with $N=100$ samples, both methods mistakenly include terms such as $[\nabla_{\mathcal{S}} u]_3$, $[\nabla_{\mathcal{S}} u]_1 [\nabla_{\mathcal{S}} u]_2$, and $u\Delta_{\mathcal{S}}u$ in the model, albeit with coefficients close to zero. However, as the number of samples increases to $N=200$, the errors in the coefficients of the terms $\Delta_{\mathcal{S}}u$ and $u$ decrease to $0.007\%$ and $0.022\%$, respectively. In this case, the normal extension method correctly excludes the remaining terms from the equation. Conversely, even with the theoretical expression of the normal vector, the term $[\nabla_{\mathcal{S}} u]_3$ remains included in the model. This discrepancy arises due to inherent errors in the node generation process, irrespective of whether the nodes are distributed on the unit sphere or any other general surfaces.


A more accurate approximation of the model can be achieved by increasing the sample size to $N=1000$ for both the normal extension method and the analytical normal vector. With the normal extension method, the errors in the coefficient of the term $\Delta_{\mathcal{S}}u$ are reduced to $1.08\times 10^{-5}\%$, and the errors of $u$ are reduced to $1.18\times 10^{-4}\%$. Meanwhile, when using the analytical normal vector, the errors are $0.007\%$ and $0.022\%$, respectively. It is worth noting that a smaller noise level in the samples can further contribute to higher-fidelity predictions. 

\begin{figure}
\centering
\includegraphics[width=3in,height=2.2in]{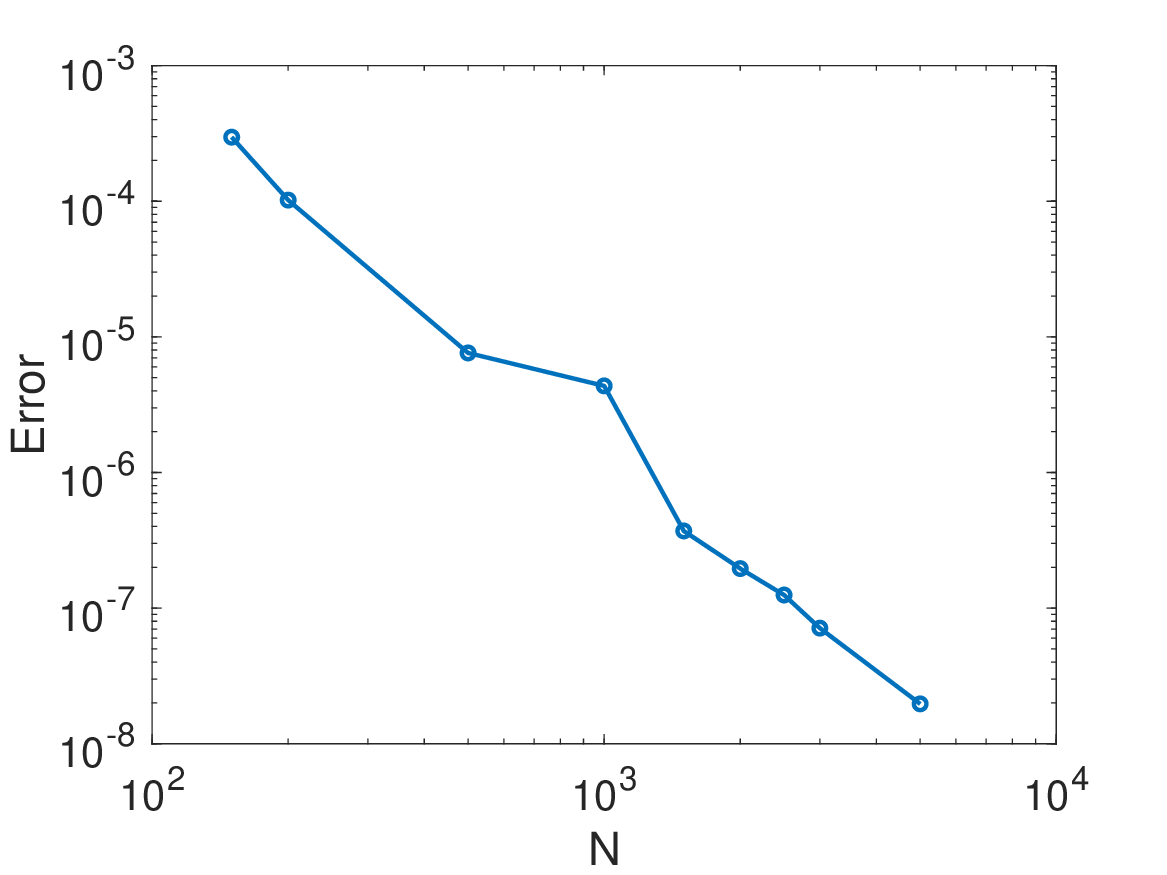}
 \caption{ The relative $L_2$ errors of solutions for solving the identified model in Example 1 on the unit sphere for different nodes sizes $N$ (in log-log scale).}
 \label{Fig:N_errEx2}
\end{figure}

\begin{figure}
	\centering
	\includegraphics[width=6in,height=5in]{
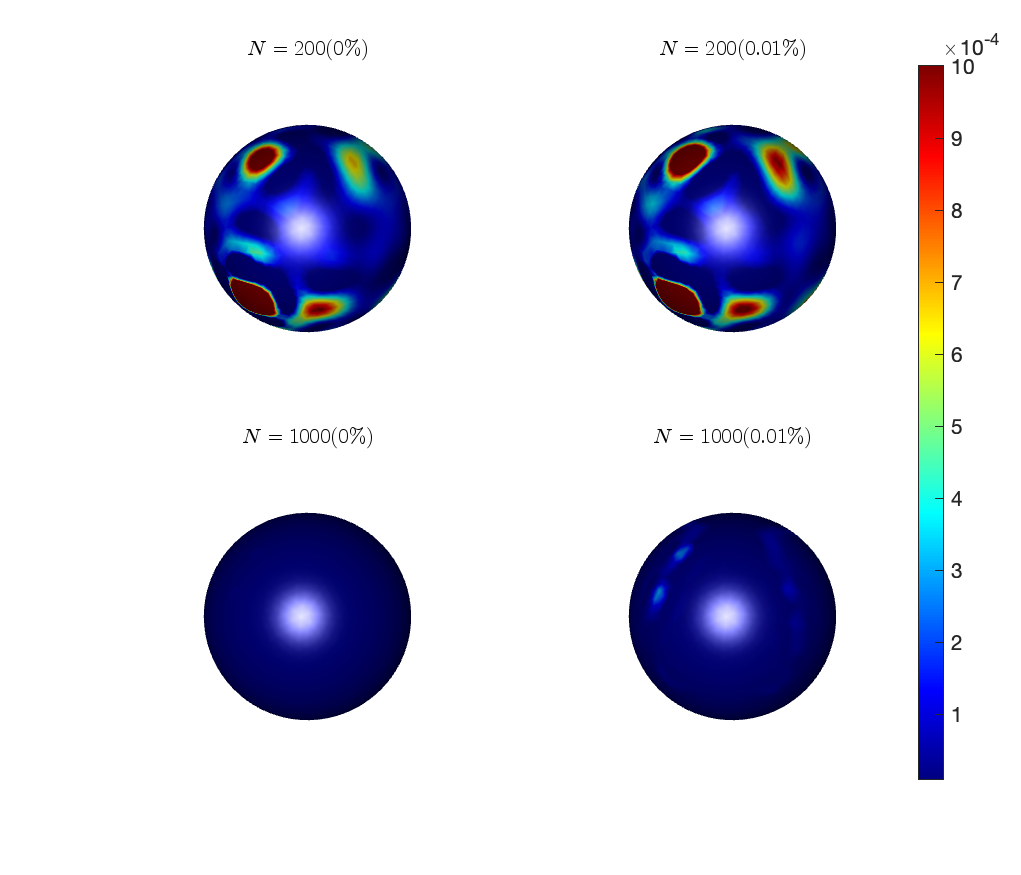 }
 \caption{ Absolute  errors of predicted solutions $u$ of the learned PDEs on the unit sphere in Example 1.}
 \label{Fig:errEx2}
\end{figure}

To assess the accuracy of the predicted solutions of the learned PDEs with the normal extension method, we provide the relative $L_2$ error in Figure \ref{Fig:N_errEx2}. This error analysis demonstrates the rapid convergence of the predicted PDE solutions. The corresponding absolute errors of the predicted solutions are presented in Figure \ref{Fig:errEx2}, which illustrate that our method can achieve an accuracy level of around $10^{-5}$ with a sample size of $1000$. These results imply that our approach performs well in accurately predicting the solution and identifying the unknown PDE simultaneously.

\begin{table}
	\footnotesize
	\centering 
	\renewcommand\arraystretch{1.5}
	\begin{tabular}{ l p{1.1cm}| l
		}
		\bottomrule[0.8pt] 
		\multicolumn{2}{c|} {\diagbox{$N$}{Noise}{Model}} &  True model  $-\Delta_{\mathcal{S}}u+u=f$ 
		\\
		\hline
		100   &  $0\%$& $-0.9992\Delta_{\mathcal{S}}u+0.9995u-0.0023 [\nabla_{\mathcal{S}} u]_3 +0.0003 [\nabla_{\mathcal{S}} u]_1 [\nabla_{\mathcal{S}} u]_2 + 0.0005 u\Delta_{\mathcal{S}}u=f$  
		\\
		\hline
		200&  $0\%$&  $-0.99985 \Delta_{\mathcal{S}}u+ 1.00043 u + 0.00024 [\nabla_{\mathcal{S}} u]_3 = f $
		\\
		\hline
		1000	&  $0\%$&  $ -0.9999975 \Delta_{\mathcal{S}}u+ 1.000071 u= f $
		\\
		&  $0.01\%$&     $-1.00007 \Delta_{\mathcal{S}}u+ 0.99978 u= f      $ 
		\\		
		\toprule[0.8pt] 
	\end{tabular}	
	\medskip
	
	\caption{
		{The identified model  in Example 1 on the unit sphere using theoretical normal vectors for different nodes sizes $N$ and different noise levels.}}
	\label{Tab:ex2_n}
\end{table}

\begin{table} 
	\footnotesize
	\centering 
		\renewcommand\arraystretch{1.5}
	\begin{tabular}{ l p{1.1cm}| l
		}
		\bottomrule[0.8pt] 
		\multicolumn{2}{c|} {\diagbox{$N$}{Noise}{Model}} &  True model  $-\Delta_{\mathcal{S}}u+u=f$ 
		\\
		\hline
	   100   &  $0\%$& $-1.0015\Delta_{\mathcal{S}}u+0.9960u-0.0014 [\nabla_{\mathcal{S}} u]_3 +0.0013 [\nabla_{\mathcal{S}} u]_1 [\nabla_{\mathcal{S}} u]_2 + 0.0003 u\Delta_{\mathcal{S}}u=f$  
		\\
		\hline
		200&  $0\%$&  $-1.0001 \Delta_{\mathcal{S}}u+ 0.9998 u= f $
		\\
		&  $0.01\%$&    $-1.00001 \Delta_{\mathcal{S}}u+ 0.99946 u= f  $ 
		\\
		\hline
		1000	&  $0\%$&  $ -1.0000001 \Delta_{\mathcal{S}}u+ 0.9999988 u= f $
		\\
		&  $0.01\%$&     $-0.99997 \Delta_{\mathcal{S}}u+ 1.00006 u= f      $ 
		\\		
			\toprule[0.8pt] 
	\end{tabular}	
	\medskip
	
	\caption{
		{The identified model  in Example 1 on the unit sphere  with normal vectors computed by normal extension method for different nodes sizes $N$ and different noise levels.}}
	\label{Tab:ex2}
\end{table}

At the end of this example, as a supplementary component to our investigation into noise robustness, we shall undertake an examination of the square-root LASSO. Then instead of Equ.(\ref{Lasso}), the sparse optimization of the coefficients will be 
\begin{equation}\label{srLasso}
\widetilde{\bm{\xi}}=\mathop{\rm argmin}\limits_{\bm{\xi}} ||\widetilde{\mathcal{L}}_{\mathcal{S}}u(X) - f(X)||_2 + \mu||\bm{\xi}||_1,
\end{equation}
and the penalty parameter $\mu$ is independent of the variance of the noise. We will utilize the recommended choice of penalty parameter in \cite{Belloni}. 
To compare with LASSO, we will test the same function  $u=10xyz+5xy+z$ satisfying the linear stationary PDE on unit sphere with different noise-level samples. Besides, the normal vectors are computed by normal extension method.  

The numerical results are presented in Table \ref{Tab:srLASSO}. For small number of sample size, for example $N=100$, even in the presence of noise level amounting to $0.1\%$ , the square-root LASSO demonstrates the advantageous property of incorporating only one spurious term $[\nabla_{\mathcal{S}} u]_3$ in the model. In contradistinction,  when operating in a noise-free environment, the classical LASSO selects a model that includes three spurious terms. However, when the sample size is increased to $N=200$,  LASSO is now able to successfully select the correct model under low noise conditions. In contrast, the square-root LASSO continues to  keep the irrelevant term  $[\nabla_{\mathcal{S}} u]_3$. As the sample size is increased even further, while the square-root LASSO is still unable to identify the true model, in the absence of noise, the coefficient of the incorrectly selected term has become quite small. 

Overall, in terms of robustness to noise, the square-root LASSO does exhibit better performance. However, when it comes to the accuracy of model learning, the standard LASSO  outperforms the square-root LASSO.

\begin{table} 
	\footnotesize
	\centering 
		\renewcommand\arraystretch{1.5}
	\begin{tabular}{ l p{1.1cm}| l
		}
		\bottomrule[0.8pt] 
		\multicolumn{2}{c|} {\diagbox{$N$}{Noise}{Model}} &  True model  $-\Delta_{\mathcal{S}}u+u=f$ 
         \\
		\hline
		100 &  $0\%$&      
          $-0.999990\Delta_{\mathcal{S}}u+1.005279u -0.002903[\nabla_{\mathcal{S}} u]_3 =f$    
         \\ &  $0.1\%$
          & $ -0.9943\Delta_{\mathcal{S}}u+1.0491u  -0.0100  [\nabla_{\mathcal{S}} u]_3 =f$    
         \\
		\hline
	   200   &  $0\%$& $-1.000042\Delta_{\mathcal{S}}u+0.999958 u -0.000099[\nabla_{\mathcal{S}} u]_3 =f$ 
		\\
		&  $0.1\%$&    $ -0.99102  \Delta_{\mathcal{S}}u+ 1.06922 u -0.00998 [\nabla_{\mathcal{S}} u]_3 = f  $ 
		\\
		\hline
		1000	&  $0\%$&  $ -0.999998 \Delta_{\mathcal{S}}u+ 1.000017  u  -0.000004[\nabla_{\mathcal{S}} u]_3= f $
		\\
		&  $0.1\%$&     $-0.9594  \Delta_{\mathcal{S}}u+ 1.3161 u -0.0518 [\nabla_{\mathcal{S}} u]_3 = f      $ 
		\\
		\hline    
		2000	&   $0\%$&  $ -0.9999994  \Delta_{\mathcal{S}}u+  1.0000050  u -0.0000009 [\nabla_{\mathcal{S}} u]_3 = f $
  \\  
  &   $0.1\%$&  $-0.86   \Delta_{\mathcal{S}}u+  2.12   u  -0.19  [\nabla_{\mathcal{S}} u]_3 = f $
  \\
			\toprule[0.8pt] 
	\end{tabular}	
	\medskip
	
	\caption{
		{The identified model computed by square-root LASSO in Example 1 on the unit sphere for different nodes sizes $N$ and different noise levels.}}
	\label{Tab:srLASSO}
\end{table}

\subsection*{}

\noindent\textbf{Example 2: Evolutionary PDE}

\vspace{1em}

Now we consider a time-dependent convective-diffusion-reaction equation on several more general surfaces, including the Torus, Brestzel2, and Cyclide, in addition to the unit sphere. The equation involves an unknown diffusive coefficient $a$, the velocity of surface motion $\bm{b}=(b_1,b_2,b_3)^\top$, a reaction coefficient $r$, and a constant $c$. The equation is given by
\begin{align*}
&\frac{\partial u}{\partial  t}=(a\Delta_{\mathcal{S}}-\bm{b}\cdot\nabla_{\mathcal{S}}+c)u+rg(u)+f, \ \  (x,y,z) \in \mathcal{S},\ \  t \in [0,T], 
\end{align*}
with the reaction term $g(u)=u^2$ and the initial condition $u|_{t=0}=u_0(x,y,z)$. 
Here $\bm{b}$ is assumed to have both the normal and tangent components, so $\bm{b}$ with any added normal component yields the same results.

The implicit expressions for the surfaces aforementioned are

$\bullet$ Torus:

$$S= (x^2+y^2+z^2+1^2-(1/3)^2)^2-4(x^2+y^2)=0 $$

$\bullet$ Cyclide:

$$ S=(x^2+y^2+z^2-1+1.9^2)^2-4(2x+\sqrt{4-1.9^2})^2-4(1.9y)^2=0$$

$\bullet$ Bretzel2:

$$ S= (x^2(1-x^2)-y^2)^2+1/2z^2-1/40(x^2+y^2+z^2)-1/40=0$$ 

In this example, the temporal snapshots of the form $\{(x_i,y_i,z_i,t_j), u^{*,j}(x_i,y_i,z_i)\}_{i, j=1}^{N, M}$ with initial data  at $t=0$ and final time $T=M\Delta t$ are given. The smoothness order $m$ of the kernels used in this example is chosen to be $4$.

\vspace{1em}


\vspace{1em}

Firstly, we consider the unit sphere and set the parameters as follows: $a=\frac{1}{2}$, $\bm{b}=\bm{0}$, $c=0$, $r=\frac{1}{8}$, and the exact solution $u=e^{x+y+z}\exp(-t)$. The time step size is $\Delta t=0.01$, and the number of time steps is $M=100$.

Table \ref{Tab:ex2_sphere} summarizes the results obtained using the normal extension method on the unit sphere, considering different node sizes and noise levels. For a node size of $N=100$, the learned model includes the term $(\Delta_{\mathcal{S}} u)^2$ with a relatively small coefficient of $1.8331\times 10^{-4}$. By increasing the node size to $N=500$, we successfully recover the model with improved accuracy. The errors in the coefficient of the term $\Delta_{\mathcal{S}}u$ reduce to $0.011\%$, and the errors in the coefficient of the reaction term $u^2$ reduce to $0.005\%$. Even when the noise level is raised to $0.01\%$, the PDE model can still be detected, albeit with slightly higher errors in the coefficients. Specifically, the errors in the coefficient of the term $\Delta_{\mathcal{S}}u$ are $0.3\%$, and the errors in the coefficient of $u^2$ are $0.7\%$.

\begin{figure}
	\centering
  \begin{overpic}[width=5.5in]{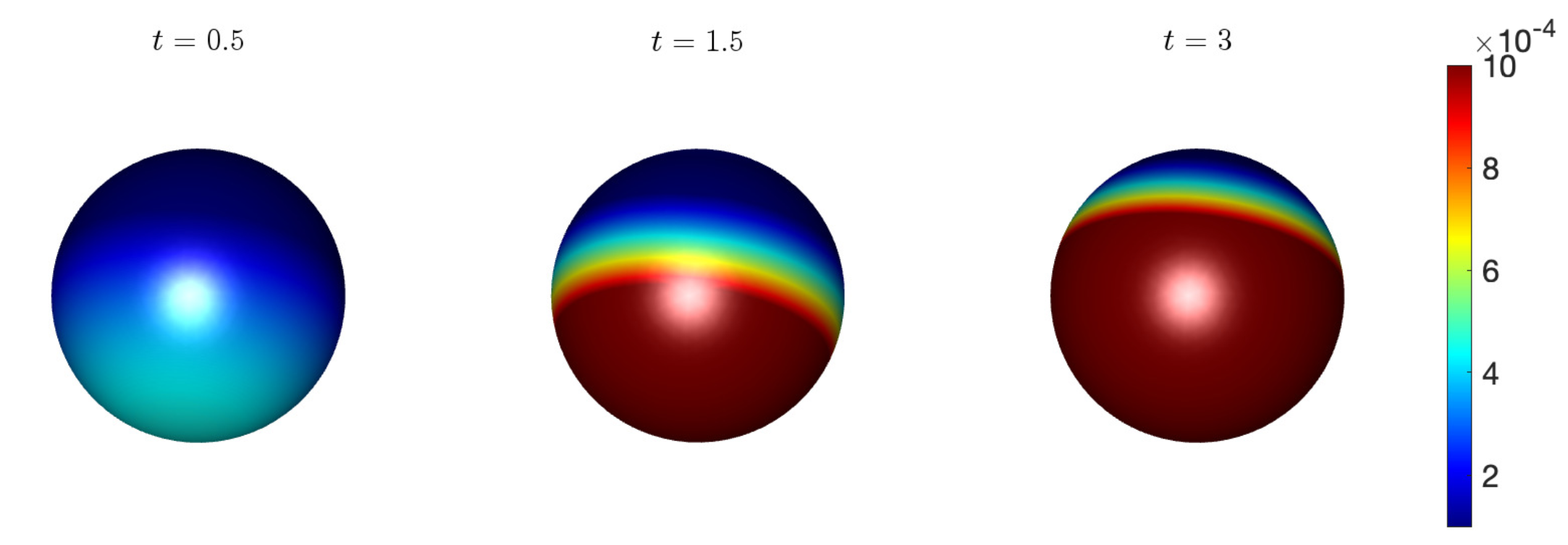}
  \tiny
    \put(0, 11) {\begin{rotate}{90}noise free\end{rotate}}
    \end{overpic}
    \begin{overpic}[width=5.5in]{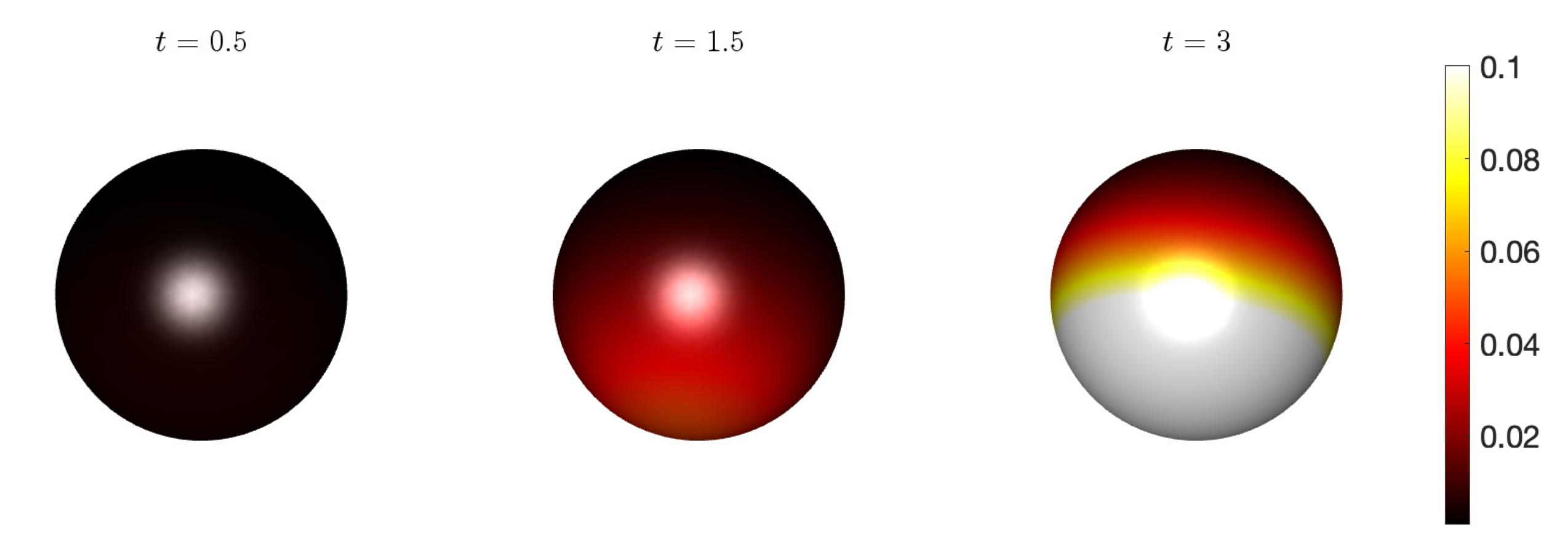}
  \tiny
    \put(0, 13) {\begin{rotate}{90}noisy \end{rotate}}
    \end{overpic}
 \caption{ Absolute  errors of predicted solutions $u$ of the learned PDEs on the unit sphere in Example 2 at different time steps with $N=500$, and $0\%$ (first row)  and $0.01\%$ (second row) noise level respectively. }\label{Fig:errEx2_sphere}
\end{figure}

Figure \ref{Fig:errEx2_sphere} provides a representation of the absolute errors in the numerical solutions $u$ obtained from the learned PDEs at different time steps. These results were obtained using a node size of $N=500$ and a noise level of $0\%$ on the unit sphere. Notably, the figure demonstrates that the method is capable of making accurate predictions even when the prediction time exceeds the final time $T$ of the available samples. This ability to accurately predict the solution over long time intervals highlights the robustness and reliability of the method.

\begin{table}
	\footnotesize
	\centering 
		\renewcommand\arraystretch{1.5}
	\begin{tabular}{ l p{1.1cm}| l 
		}
 \bottomrule[0.8pt] 
		\multicolumn{2}{c|} {\diagbox{N}{noise}{model}}&  True model  $\frac{\partial u}{\partial  t}=0.5\Delta_{\mathcal{S}}u+0.125u^2+f$
		   \\
		   \hline
		   100 &  $0\%$&  $\frac{\partial u}{\partial  t}=0.50001\Delta_{\mathcal{S}}u+0.12499u^2+1.8331\times 10^{-4}(\Delta_{\mathcal{S}} u)^2+f$
		   \\
		   \hline
	          500 &  $0\%$  &  $\frac{\partial u}{\partial  t}=0.50005\Delta_{\mathcal{S}}u+0.12499u^2+f$
		  	\\
		  &  $0.01\%$&  $\frac{\partial u}{\partial  t}=0.49859\Delta_{\mathcal{S}}u+0.12413u^2+f$
		  \\
		  \hline
		     1000 &  $0\%$  &  $\frac{\partial u}{\partial  t}=0.50006\Delta_{\mathcal{S}}u+0.12500u^2+f$
		  \\
		  &  $0.01\%$&  $\frac{\partial u}{\partial  t}=0.49356\Delta_{\mathcal{S}}u+0.12033u^2+f$
		  \\
	 \toprule[0.8pt] 
		  
	\end{tabular}	
	\medskip
	\caption{\small
		{The identified model in learning convective-diffusion-reaction equation on the unit sphere for different nodes sizes $N$ with different noise levels. ($\Delta t = 0.01$)}}
	\label{Tab:ex2_sphere}
\end{table}

\vspace{1em}



In the case of general surfaces, we test another parameter set: $a=1$, $\bm{b}=\bm{0}$, $c=0$, and $r=1$. The exact solution is assumed to be $u=\sin(x)\sin(y)\sin(z)\sin(t)$. The time step size is $\Delta t = 0.01$, and the number of time steps is $M=100$ as before. 

The complexity of the surfaces impacts the number of nodes required to accurately recover the underlying models. Due to the distinct characteristics of each feature in the surface PDE, the regularization parameters in the algorithms can be chosen to be quite small or even zero. In fact, for $L_1$ regularization, coefficients below a certain tolerance can be simply deleted.

On the torus, with a node size of $N=3968$, we obtain the model with errors in the parameter $a$ being $0.0098\%$, and errors in $r$ being $2.13\times 10^{-4}\%$. Increasing the number of samples to $N=7520$ has a minimal effect on improving the accuracy of the model recovery. Similar results are observed with nodes on the one-holed Cyclide.

\begin{figure} [htp]
	\centering
	
		\includegraphics[width=6in,height=6in]{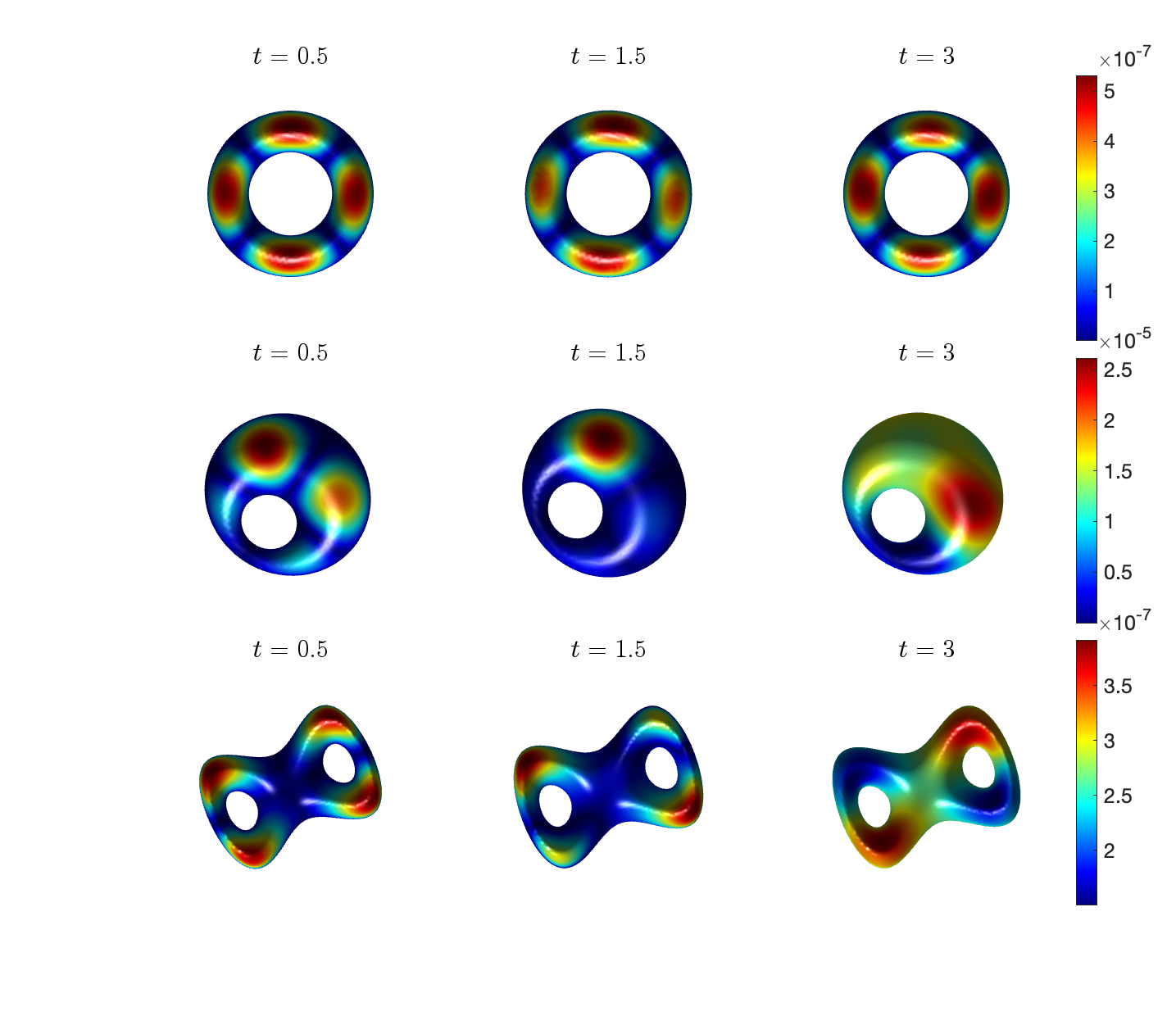}

	\caption{ The  absolute  errors  of predicted solutions from the learned PDE in Example 2 on general surfaces with samples of $0\%$ noise level. (top) Torus  with $N=3968$. (middle) Cyclide with $N=3662$. (bottom) Bretzel2 with $N=9546$.
}\label{Fig:errEx2_surfaces} 
\end{figure}

However, in the case of the $2$-holed Bretzel2, even more nodes ($N=7270$)  are  required to accurately reconstruct the underlying PDE.  Here, the errors in the parameters $a$ and $r$ are reduced to $0.0093\%$ and $0.0633\%$, respectively. More details of the specific learned models can be found in Table \ref{Tab:ex3}.
Figure \ref{Fig:errEx2_surfaces} shows the absolute errors of the predicted solution obtained from the learned PDEs with nodes on the general surfaces. It can be observed that accurate long-time predictions can also be achieved.

\begin{table} 
 \footnotesize
	\centering 
	\renewcommand\arraystretch{1.5}
	\begin{tabular}{ l l | l 
		}
		\hline
		\multicolumn{2}{l|} {\diagbox[width=9em]{Surface}{$N$}{Model}} &  True model  $\frac{\partial u}{\partial  t}=\Delta_{\mathcal{S}}u+u^2+f$
		\\
		\hline
		 Torus & $3968$& $\frac{\partial u}{\partial  t}=1.000098\Delta_{\mathcal{S}}u+0.999993u^2+f$
		\\
		 &  $7520$&  $\frac{\partial u}{\partial  t}=1.000098\Delta_{\mathcal{S}}u+1.0000025u^2+f$\\
		\hline
	    Cyclide& $3662$    &  $\frac{\partial u}{\partial  t}=1.000093\Delta_{\mathcal{S}}u+1.00000002u^2+f$
		\\
	    	&$7030$ &  $\frac{\partial u}{\partial  t}=1.000098\Delta_{\mathcal{S}}u+0.9999998u^2+f$
	    	\\
		\hline
		Bretzel2 & $7270$     &  $\frac{\partial u}{\partial  t}=1.000093\Delta_{\mathcal{S}}u+1.000633u^2+f$ \\
		 &$13014$&  $\frac{\partial u}{\partial  t}=1.000097\Delta_{\mathcal{S}}u+0.999960u^2+f$ 
		\\
		\hline
	\end{tabular}	
	\medskip
	
	\caption{
		{The identified model in learning convective-diffusion-reaction equation on general surfaces  for different nodes sizes $N$. ($\Delta t = 0.01$)}}
	\label{Tab:ex3}
\end{table}


\vspace{1em}

\noindent\textbf{Example 3: Eikonal equation}

\vspace{0.8em}

In this example, we will demonstrate the numerical approximation of the geometric distance using the large $p$-Laplacian, which allows us to learn the eikonal equation \cite{plap}.

The $p$-Laplacian operator on surface $\mathcal{S}$ is defined as  
$$
\Delta_{\mathcal{S}}^p u = \nabla_{\mathcal{S}}\cdotp(\left|\nabla_{\mathcal{S}}u\right|^{p-2}\nabla_{\mathcal{S}}u). 
$$
The parameter $p$ determines the degree of nonlinearity, with larger values of $p$ placing more emphasis on capturing sharp features and geometric properties. We will focus on the range $p \in [2, +\infty)$.

\begin{figure}
\vspace{-1cm}
\centering
\begin{minipage}[t]{0.48\textwidth}
\centering
\includegraphics[width=3in]{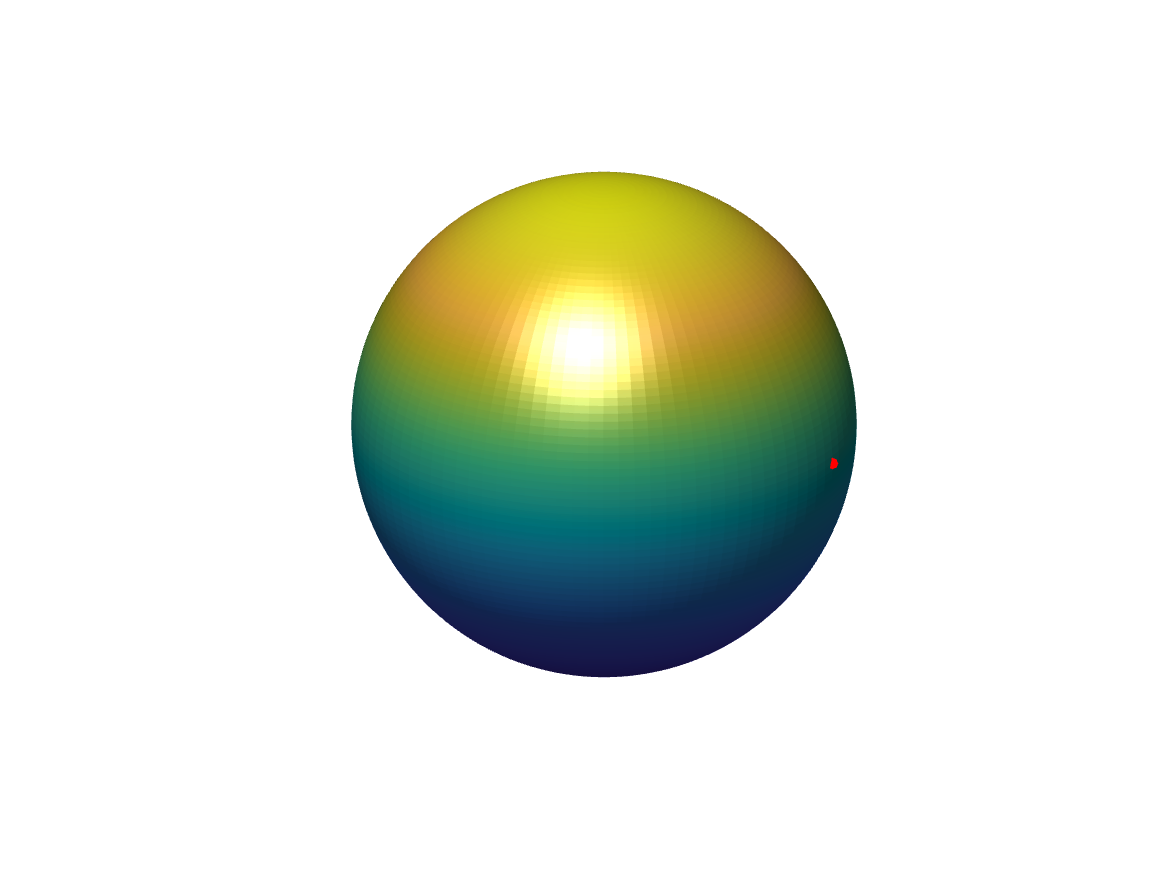}
\caption*{}
\end{minipage}
\begin{minipage}[t]{0.48\textwidth}
\centering
\includegraphics[width=3in]{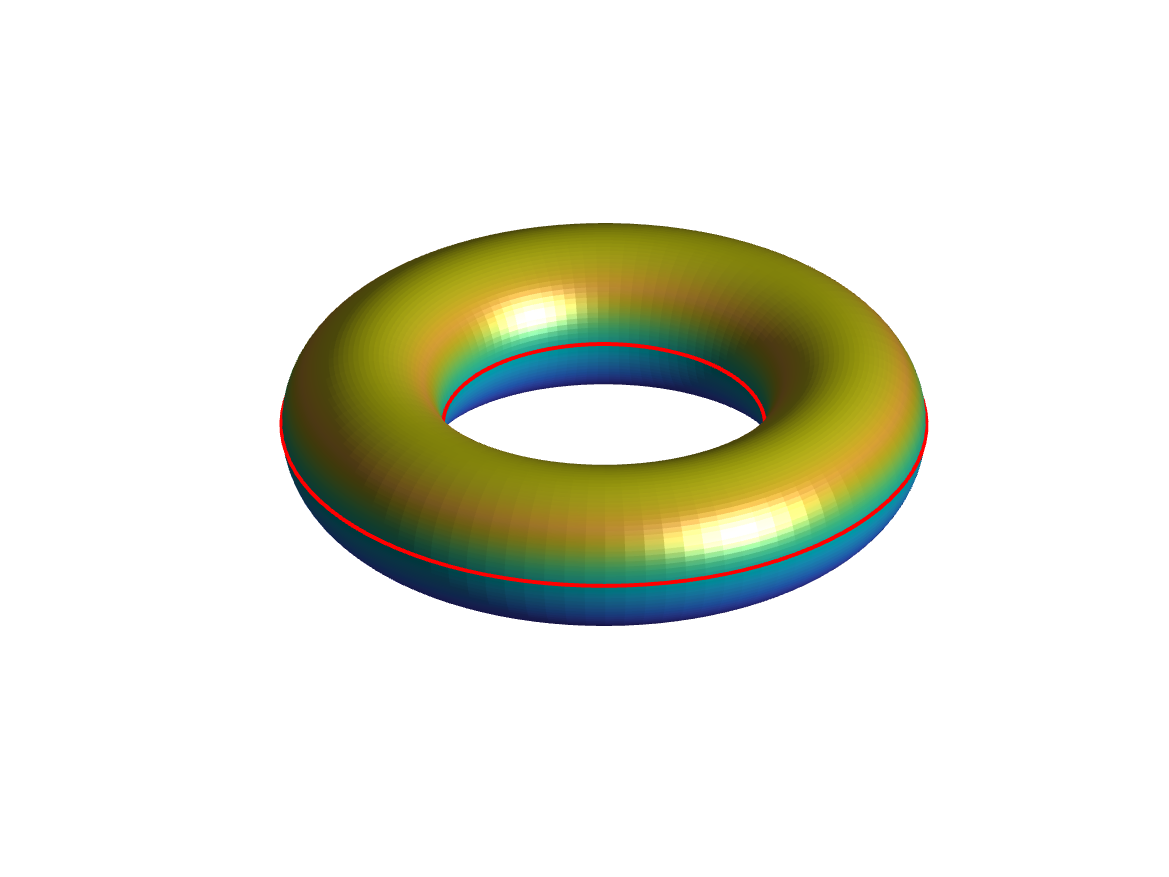}
\caption*{}
\end{minipage}
\vspace{-1cm}
 \caption{(left) The unit sphere with the source point at $\Sigma =\bm{x}_s=(0,0,1)$ (indicated in red) and (right) the torus with source points being the two closed circles (indicated in red).}
    \label{fig:sphere_torus_EEp}
\end{figure}

The eikonal equation on the surface in a special case is given by 
\begin{align}
     & \left|\nabla_{\mathcal{S}}u\right|= 1, \ \  \bm{x} \in \mathcal{S} \backslash \Sigma  \\  
     & u(\bm{x}_s) = 0,  \ \ \bm{x}_s \in \Sigma   \notag
\end{align}
where $\Sigma$ is the set of sources, and $u(\bm{x})$ represents the geodesic distance from the target point $\bm{x} \in \mathcal{S} \backslash \Sigma$ to $\Sigma$.
Specific choices of source points on the sphere and torus can be seen in Figure \ref{fig:sphere_torus_EEp}. The eikonal equation has various applications, including image regularization \cite{liu} and modeling constrained motions of curves or surfaces.

When dealing with the additional condition imposed on the source points, we can employ a two-step learning process to approximate the eikonal equation using a library of tailored $p$-Laplacian operators for surfaces.
In the first step of the learning process, we aim to learn an initial approximation of the eikonal equation without considering the constraint at the source points. This step involves training a model using the available samples $\{\bm{x}_i,u^*_i\}_{i=1}^{N} $, where $u_i^*$ represents the geodesic distance from $\bm{x}_i$ to the set of sources. The objective is to minimize the discrepancy between the predicted distance function and the given geodesic distances for the non-source points.


First, we construct the feature library $\Lambda(X)$, which includes the linear polynomial of $u$ and the $p$-Laplacian with different $p$ values evaluated at the nodes $X$, 
\begin{equation*}\label{library_EE}
\Lambda(X)= \Big[u|_X, \Delta_{\mathcal{S}}^{p_1} u|_X, \dots, \Delta_{\mathcal{S}}^{p_m} u|_X \Big], 
\end{equation*}
where $\{p_1, p_2, \dots, p_m\}$ are the candidate $p$ values.  
To approximate the feature library, we use a similar approach as in (\ref{gradient1}) and (\ref{laplace}). The model can be written as 
$$
\widetilde{\Lambda}(X) \bm{\xi}= f(X),
$$
with $f(X)=\bm{1}$  being the unit column vector.   
By solving the sparse optimization problem (\ref{Lasso}), we obtain the coefficients $\widetilde{\bm{\xi}}$.
Once we have an initial approximation, we calculate the model residuals
$$\bm{r}(X)=\widetilde{\Lambda}(X) \widetilde{\bm{\xi}}-f(X).$$
with $\bm{r}(X)=[r(\bm{x}_1),r(\bm{x}_2),\dots,r(\bm{x}_N)]^{\top}$.
In the second step, we fit the residuals obtained from the first step to identify the source points. We model the residuals as a sum of Gaussian kernel functions centered at the potential source points 
 $$
 r(\cdot)=\sum_{i=1}^{N}\eta_i \Psi(\cdot-x_i) 
 $$
where $\Psi(\cdot-x_i)=e^{-\frac{||\cdot-x_i||}{\sigma^2}}$ is the Gaussian kernel function  restricted on surfaces with the shape parameter $\sigma^2$, which can simulate the delta function located at the potential source point, and $\bm{\eta}=[\eta_1,\eta_2,\dots,\eta_N]^{\top}$ is the unknown coefficients of the source terms. To determine the sparse source terms, we solve the following optimization problem 
\begin{equation*}\label{Lasso2}
\widetilde{\bm{\eta}}=\mathop{\rm argmin}\limits_{\bm{\eta}} ||\Psi(X,X)\bm{\eta} - \bm{r}(X)||_2^2 + \mu||\bm{\eta}||_1.
\end{equation*}

\vspace{1em}


In the first experiment, we test the nodes on the unit circle with different one-point sources. We select potential values of $p$ from ${2, 5, 50, 100, 200, \dots, 1000}$. For $N=100$, we choose the source term to be $\bm{x}_s = \bm{x}_{20}$, and for $N=1000$, we choose $\bm{x}_s = \bm{x}_{1}$.

The results of the experiments are summarized in Table \ref{Tab:exEcirle}. We find that except for the $u$ term and a low $p$ value ($p=2$) term, the high $p$ value ($p=1000$) is consistently selected in the model. This suggests that larger $p$ values are preferred for accurately learning the eikonal equation. The learned right-hand side term in the equation includes $\bm{x}_{16}$ (near the source) for $N=100$ and $\bm{x}_{1}$ (exactly the source) for $N=1000$. These results demonstrate that our algorithm can learn the eikonal equation using a large $p$-Laplacian and simultaneously identify the source point.


\begin{table}
	\footnotesize
 \centering 
	\renewcommand\arraystretch{1.5}
	\begin{tabular}{  l l|l 
		}
	\bottomrule[0.8pt]
			\multicolumn{2}{ c|} {\diagbox[width=9em]{$N$}{Source}{Model}}   & True model 
    $\left|\nabla_{\mathcal{S}}u\right|= 1, \ \  \bm{x} \in \mathcal{S} \backslash \{\bm{x}_s\} $

  \\
		\hline
		 $100$   & $  \bm{x}_s=\bm{x}_{20}$ &  $2.6678 \Delta_{\mathcal{S}}^{1000} u +2.6745 \Delta_{\mathcal{S}}u + 0.9585u -1 =-0.0110e^{-||\bm{x}-\bm{x}_{16}||^2} $
		\\
		 $1000$  &  $  \bm{x}_s=\bm{x}_{1} $ &  $0.0472 \Delta_{\mathcal{S}}^{1000} u  -0.0034  \Delta_{\mathcal{S}}u + 0.9245u -1 =-0.1277e^{-||\bm{x}-\bm{x}_1||^2}$    
		\\ 

			\toprule[0.8pt] 
	\end{tabular}	
	\medskip
	
	\caption{
 The identified model in Example 3 on the unit circle with different numbers of samples and different source points. Here $\bm{x}_{1}=(1,0)$, $\bm{x}_{16}=(0.5801,0.8146)$, and $\bm{x}_{20}=(0.3569,0.9341)$. 
  }
	\label{Tab:exEcirle}
\end{table}

\vspace{1em}


Next, we consider the surface $\mathcal{S}$ to be the unit sphere and impose a boundary condition at a single point source $\bm{x}_s$. We conduct experiments with different numbers of samples ($N=1000$ and $N=2000$), and the potential values of $p$ were selected from ${2, 5, 50, 100, 200, \dots, 1000}$. The results of the experiments are presented in Table \ref{Tab:exEsphere}. To deal with the ill-conditioned problem arising from large $p$ values, we increased the number of samples to obtain a more reliable approximation model. For $N=1000$, the right-hand term in the model was $\bm{x}_1$, which corresponds to the exact source point. However, for $N=2000$, the right-hand term was $\bm{x}_{85}$, representing a point near the source. It is worth noting that regardless of the specific source point and the number of samples, the model consistently identified the same $p$-Laplacian terms. This indicates that the choice of a high $p$ value ($p=1000$) is robust and independent of the source point and the number of samples.

\begin{table}
	\footnotesize
 \centering 
	\renewcommand\arraystretch{1.5}
	\begin{tabular}{  l l|l 
		}
	\bottomrule[0.8pt]
		\multicolumn{2}{ c|} {\diagbox[width=9em]{$N$}{Source}{Model}}  & 
  True model 
    $\left|\nabla_{\mathcal{S}}u\right|= 1, \ \  \bm{x} \in \mathcal{S} \backslash \{\bm{x}_s\} $
  \\
		\hline
		 $1000$   & $  \bm{x}_s=\bm{x}_{1}$ &  $-0.0155 \Delta_{\mathcal{S}}^{1000} u +0.2616 \Delta_{\mathcal{S}}u + 1.0302u -1 = 0.0001e^{-||\bm{x}-\bm{x}_{1}||^2} $
		\\
		 $2000 $  &  $  \bm{x}_s=\bm{x}_{20} $ &  $-0.0018 \Delta_{\mathcal{S}}^{1000} u  +0.2854  \Delta_{\mathcal{S}}u + 1.0344u -1 = -0.0009e^{-||\bm{x}-\bm{x}_{85}||^2}$    
		\\ 

			\toprule[0.8pt] 
	\end{tabular}	
	\medskip
	\caption{
		The identified model in Example 3 on the unit sphere with different numbers of samples and different source points.
Here $\bm{x}_{1}=(0,0,1)$, $\bm{x}_{20}=(-0.1355,0.0782,0.9877)$, and $\bm{x}_{85}=(-0.2676,0.1545,0.9511)$. 
}
	\label{Tab:exEsphere}
\end{table}

\begin{table}
	\footnotesize
 \centering 
	\renewcommand\arraystretch{1.5}
	\begin{tabular}{   l|l 
		}
	\bottomrule[0.8pt]
		 \diagbox[width=5.5em]{$N$}{Model} & 
  True model 
    $\left|\nabla_{\mathcal{S}}u\right|= 1, \ \  \bm{x} \in \mathcal{S} \backslash \Sigma $
  \\
		\hline
		 $3968$   
   &  $ -0.0005\Delta_{\mathcal{S}}^{100} u + 0.0087\Delta_{\mathcal{S}}^{15}u + 3.0803u -1 = -0.0011e^{-||\bm{x}-\bm{x}_{1807}||^2} $
		\\
		 $10814 $  & 
   $-0.0012\Delta_{\mathcal{S}}^{1000} u  + 0.0043\Delta_{\mathcal{S}}^{400}u + 3.0100u -1 = -0.0008e^{-||\bm{x}-\bm{x}_{5281}||^2}$    
		\\  

			\toprule[0.8pt] 
	\end{tabular}	
	\medskip
	\caption{
		{The identified model in Example 3 on the torus with different numbers of samples. Here $\bm{x}_{1807}=(0.4472,-0.4944,0.0005)$, and $\bm{x}_{5281}=(-0.2588,0.6144,0.0005)$. }}
	\label{Tab:exEtorus}
\end{table}

\vspace{1em}


Finally, we test the geodesic distance to $P$ containing two closed curves on a torus, see the right side of Figure \ref{fig:sphere_torus_EEp}. The torus has an outer radius of $1$ and an inner radius of $\frac{1}{3}$, and $P$ contains the inner and outer circles lying in the plane $z=0$. The experiments are conducted with $N=3968$ and $N=10814$ samples, and the potential values of $p$ are selected from ${2, 5, 10, 15, \dots, 100}$ and ${2, 5, 50, 100, \dots, 1000}$, respectively.
The results of the experiments are shown in Table \ref{Tab:exEtorus}. Similar to the previous example, when dealing with more complex surfaces, a larger number of samples is required to effectively learn the model and handle the ill-conditioned problem caused by large values of $p$. In both cases, the models consistently include the largest value of $p$ among the selected potential values. The right-hand term in the model is $\bm{x}_{1807}$ for $N=3968$ and $\bm{x}_{5281}$ for $N=10814$, which are both close to the source circles of the torus.

\vspace{1em}

\noindent\textbf{Example 4: Biharmonic equation}

\vspace{0.8em}

The above examples only include the surface PDE operators $\mathcal{L}_{\mathcal{S}}u := \mathcal{L}_{\mathcal{S}}(u,\nabla_{\mathcal{S}}u,\Delta_{\mathcal{S}}u)$ up to second-order. 
Actually, the aforementioned methodology can be generalized to accommodate higher-order surface PDEs. As an illustrative case, here we shall examine the biharmonic equation on unit sphere. 

The biharmonic equation arises in the modeling of more complex phenomena in solid mechanics and fluid mechanics. For instance, 
 boundary value problems for the biharmonic equation were used to model radar imaging with broad-band and low-frequency waves \cite{Andersson}, and 
the stream function of incompressible Stokes flow in two-dimensional space is the solution of a biharmonic equation \cite{Lai}. 

Let $u_l^m(\theta, \phi)$ being the spherical harmonic function of degree $l$ and order $m$ ($-l \leq m \leq l$) with  $\theta$ and $\phi$ representing colatitude and longitude, respectively. $u_l^m(\theta, \phi)$ can be defined by a product of trigonometric functions and associated Legendre polynomials. 
Then it can be verified that $u_l^m$ satisfies the following  biharmonic equation on sphere

\begin{equation}\label{biharmonic}
    \Delta^2_{\mathcal{S}}u = l^2(1+l)^2 u
\end{equation}
with $\Delta^2_{\mathcal{S}}= \Delta_{\mathcal{S}}(\Delta_{\mathcal{S}})$.  With the samples $\{X, u^{*}(X)\}$ and $X=\{(x_i,y_i,z_i)\}_{i=1}^{N}$ on unit sphere, our aim is then to identify this surface PDE from the candidate equation
$$
\mathcal{L}_{\mathcal{S}}u:=\mathcal{L}_{\mathcal{S}}(u,\nabla_{\mathcal{S}}u,\Delta_{\mathcal{S}}u, \nabla_{\mathcal{S}}(\Delta_{\mathcal{S}}u),
\Delta^2_{\mathcal{S}}u) = u. 
$$

To get the discrete approximations to the higher order surface differential operators by using radial basis function approximation methods, 
the technique employed for the function itself can also be applied to Laplaca-Beltrami operator. Specifically,  with the ‘samples’ $\widetilde{\Delta}_{\mathcal{S}} u(X)$,  we replace $u$ by $\widetilde{\Delta}_{\mathcal{S}}u$ from Equ.(\ref{interpolant}) to 
Equ.(\ref{laplace2}). Then we will obtain the interplolant $\widetilde{\widetilde{\Delta}_{\mathcal{S}} u}$, the estimate of the surface gradient of the Laplaca-Beltrami operator $\widetilde{\nabla}_{\mathcal{S}} (\widetilde{\Delta}_{\mathcal{S}} u)$, and in the end $\widetilde{\Delta}_{\mathcal{S}} (\widetilde{\Delta}_{\mathcal{S}} u)$.
Accordingly, the feature matrix $\Lambda(X)$ up to second-order degree polynomial can be constructed with 
{\small
$$
\bm{\lambda}(\bm{x})=(u(\bm{x}),
[\nabla_{\mathcal{S}} u]_1(\bm{x}),\ldots,
[\nabla_{\mathcal{S}} u]_d(\bm{x}),\Delta_{\S}u(\bm{x}),
[\nabla_{\mathcal{S}} (\Delta_{\S}u)]_1(\bm{x}),\ldots, 
[\nabla_{\mathcal{S}} (\Delta_{\S}u)]_d(\bm{x}), \Delta^2_{\mathcal{S}} u).
$$
}

Next we will choose the solution of the surface PDE to be $u:=u_1^0=\cos(\theta)$ to validate the experiment with the noise free samples, and the results are shown in Table \ref{Tab:ex4}. 
In this case, the features in the equation will be  more than $50$ terms, 
even though the biharmonic equation can be learned even by a quite small number of sample size $N=50$. Besides, it only results in a error $0.006\%$ of the coefficient. 
When the sample size reaches $1000$, we will obtain a true model with a negligible error. 
However, it is well known that obtaining precise approximations of high-order derivatives is an exceedingly arduous task. Consequently, procuring accurate model estimates for sample data exhibiting error is not a trivial endeavor. 

\begin{table} [H]
  \footnotesize
	\centering 
		\renewcommand\arraystretch{1.5}
	\begin{tabular}{ l| l
		}
		\bottomrule[0.8pt] 
	    \diagbox[width=5.5em]{$N$}{Model}&  True model  $0.25\Delta^2_{\mathcal{S}}u = u$
         \\
		\hline
	50  & $0.249985\Delta^2_{\mathcal{S}}u = u$  
		\\
	  100  & $ 0.2499923\Delta^2_{\mathcal{S}}u = u$  
		\\
		200&   $ 0.2499978\Delta^2_{\mathcal{S}}u = u$
		\\
		1000 &     $ 0.249999998\Delta^2_{\mathcal{S}}u = u $
			\\
   \toprule[0.8pt] 
	\end{tabular}	
	\medskip
	
	\caption{
		{The identified model in Example 4 on the unit sphere  with normal vectors computed by normal extension method for different nodes sizes $N$ and noise free samples.}}
	\label{Tab:ex4}
\end{table}

\section{Conclusions}

In this paper, we have introduced the physics-informed sparse optimization (PIS) approach to uncover hidden PDE models on surfaces. The method combines the $L_2$ physics-informed model loss with an additional $L_1$ regularization penalty term in the loss function.  The inclusion of the $L_1$ regularization aids in the identification of specific physical terms within the surface PDE. The corresponding loss functions are formulated in equations (\ref{Lasso}) and (\ref{t_Lasso}). To learn both stationary and evolutionary PDE models, Algorithm 1 and Algorithm 2 are presented, respectively. 

The effectiveness of the proposed approach is validated through numerical experiments conducted in Section 4. 
These experiments involve predicting solutions and simultaneously identifying unknown PDEs on spheres as well as various curved surfaces. The results demonstrate the good performance and capability of the approach in handling different surface geometries and accurately recovering the underlying PDE models. 

In our future work, we have several plans to enhance and extend the methodology presented in this paper. 
Firstly, we intend to apply the physics-informed sparse optimization approach to surface PDEs in a wider range of application domains, such as biology and image processing. 
One important aspect we will focus on is improving the stability of the algorithm and its ability to handle noisy data in practical scenarios. 
Additionally, we plan to extend the scope of our approach to more complex surfaces and even high-dimensional manifolds. This expansion will enable us to tackle challenging problems posed by the big data era, where datasets often exhibit intricate geometric structures.


\appendix




\begin{thebibliography}{99}

	
\bibitem{Adalsteinsson}
	D. Adalsteinsson, J.A. Sethian. Transport and diffusion of material quantities on propagating interfaces via level set methods.  J. Comput. Phys. 185(1), 2003, 271-288.

\bibitem{SRLasso} 
B. Adcock, A.Y. Bao, S. Brugiapaglia.  Correcting for unknown errors in sparse high-dimensional function approximation. Numer. Math. (Heidelb.) 142, 2019, 667-711.


\bibitem{Andersson}
L.E. Andersson, T. Elfving, G.H. Golub. Solution of biharmonic equations with application to radar
imaging. J. Comput. Appl. Math. 94(2), 1998, 153–180.

 
	\bibitem{Auer}
     S. Auer, R. Westermann.  A semi-Lagrangian closest point method for deforming surfaces. Computer Graphics Forum, Wiley Online Library.
	32, 2013, 207-214.

\bibitem{Belloni}
A. Belloni, V. Chernozhukov, L. Wang. Square-root lasso: pivotal recovery of sparse signals via conic programming. Biometrika. 98(4), 2011, 791–806.
 
	\bibitem{Biddle}
	H. Biddle, I. von Glehn, C. B. Macdonald, T. M{\rm $\ddot{a}$}rz. A volume-based method for denoising on curved surfaces. 
	2013 IEEE International Conference on Image Processing, IEEE.  2013, 529-533. 
	
	\bibitem{Brunton}
	S.L. Brunton, J.L. Proctor, J.N. Kutz. Discovering governing equations from data by sparse identification of nonlinear dynamical systems. 
	Proc. Natl. Acad. Sci. 113(15),  2016, 3932-3937.
	
	
	
	\bibitem{chen}
	M. Chen, L. Ling. Extrinsic meshless collocation methods for PDEs on manifolds.  SIAM J. Numer. Anal. 58(2), 2020, 988-1007.
	
	
	\bibitem{cheung}
	K.C. Cheung, L. Ling. A kernel based embedding method and convergence analysis for surface PDEs. SIAM J. Sci. Comput. 40(1), 2018, 
	 A266-A287.
	 
	 
	
	 \bibitem{Charles}
	 C.M. Elliott, B. Stinner. 
	 Modeling and computation of two phase geometric biomembranes using surface finite elements. 
	J. Comput. Phys. 
	 229(18), 
	 2010,
	 6585-6612. 
	 
	 \bibitem{Friedman}
	 J. Friedman, T. Hastie , R. Tibshirani. Regularization paths for generalized linear models via coordinate decent. J. Stat. Softw. 33(1), 2010, 1–22.
	 
	 \bibitem{Floater}
	M.S. Floater, K. Hormann. Surface Parameterization: a Tutorial and Survey. 
	Advances in Multiresolution for Geometric Modelling. Mathematics and Visualization. Springer, Berlin, Heidelberg. 2005.
	 
	 \bibitem{Fuselier}
	 E. Fuselier, G.B. Wright. Scattered data interpolant on embedded submanifolds with
	 restricted positive definite kernels: Sobolev error estimates. SIAM J. Numer. Anal. 50, 
	 2012, 1753-1776. 

      \bibitem{Fuselier13}
	 E. Fuselier, G.B. Wright, 
A high-order kernel method for diffusion and reaction-diffusion equations on surfaces.  J. Sci. Comput. 56, 2013, 535–565.
  
	 \bibitem{GaoZhen}
	 Z. Gao, Y.F. Lin, X. Sun, X.Y. Zeng.
	 A reduced order method for nonlinear parameterized partial differential equations using dynamic mode decomposition coupled with k-nearest-neighbors regression. 
	  J. Comput. Phys.
	  452, 2022, 110907.
	 
	 \bibitem{Gia}
	 Q.T. Le Gia, I.H. Sloan, H. Wendland. 
	 Multiscale analysis in Sobolev spaces on the sphere. SIAM J. Numer. Anal. 48:6, 2010, 2065-2090.
\bibitem{liu}
  J. Liu, S. Leung. A splitting algorithm for image segmentation on manifolds represented by the grid
based particle method. J. Sci. Comput. 56(2), 2013, 243–266. 
	 
	 \bibitem{Hildebrandt}
	 F. Hildebrandt, M. Trabs.
	 Nonparametric calibration for stochastic reaction–diffusion equations based on discrete observations.
	 Stoch. Process. Their Appl. 
	 162,
	 2023,
171-217.

 \bibitem{Lai}
M-C. Lai, H-C. Liu. 
Fast direct solver for the biharmonic equation on a disk and its application to incompressible flows. 
Appl. Math. Comput.
164(3), 
2005,
679-695. 

  
	   \bibitem{Lei}
	  G.H. Lei, Z. Lei, L. Shi, C.Y. Zeng,  D.-X. Zhou.   
	 Solving PDEs on spheres with physics-informed convolutional
	 neural networks. 2023. arXiv preprint, https://arxiv.org/pdf/2308.09605.pdf 
	 
	  \bibitem{Dong}
	 Z.C. Long, Y.P. Lu, X.Z. Ma, B. Dong. PDE-Net: Learning PDEs from data. 
	 PMLR 80, 2018, 3208-3216.

	 
	 \bibitem{Osborne95}
	 M.R. Osborne, G.K. Smith. A modified Prony algorithm for exponential function fitting. SIAM J. Sci. Comut. 16(1), 1995, 119-138. 
	 
	  \bibitem{mesh}
	 P.-O. Persson, G. Strang. A simple mesh generator in MATLAB. 
	 SIAM Review, 46(2), 2004, 329-345. 
	 
	 \bibitem{Ruuth}
	 A. Petras, L. Ling, S.J. Ruuth.  
	 An RBF-FD closest   point method for solving PDEs on surfaces. J. Comput. Phys.
	 370(1), 2018, 43-57.
	 
	\bibitem{OGM}
	 C. Piret. The orthogonal gradients method: A radial basis functions method for solving partial differential equations on arbitrary surfaces.  J. Comput. Phys. 
	 231(14), 2012, 4662-4675.

\bibitem{plap}
H. Potgieter,  R.C. Fetecau,  and S.J. Ruuth.  Geodesic distance approximation using a surface finite element method for the p-Laplacian. 2023. Available at SSRN: https://ssrn.com/abstract=4611037 or http://dx.doi.org/10.2139/ssrn.4611037
  
\bibitem{PINN}
M. Raissi, P. Perdikaris, G.E. Karniadakis. 
	 Physics-informed neural networks: A deep learning framework for solving forward and inverse problems involving nonlinear partial differential equations. 
	 J. Comput. Phys. 
	 378(1),
      2019,
	  686-707.
	 
	 
	 \bibitem{Ruuth1}
	 S.J. Ruuth. 
	 Implicit-explicit methods for reaction-diffusion problems in pattern formation. J. Math. Biol. 34, 1995, 148-176.
	 
	 
	 
	 \bibitem{Schaeffer}
	 H.  Schaeffer. Learning partial differential equations via data discovery and sparse optimization. Proc. Math. Phys. Eng. Sci. 473(2197), 2017, 20160446.
	 


  
	 
	 \bibitem{Schmidt}
	 M. Schmidt. Least squares optimization with L1-norm regularization. 2005. CS542B Project Report.
	 
	  \bibitem{Science}
	  M. Schmidt, H. Lipson. Distilling Free-Form Natural Laws from Experimental Data. Science 324, 2009, 81-85.


  \bibitem{Shekarpaz}
	 S. Shekarpaz, F. Zeng and G. Karniadakis. 
Splitting physics-informed neural networks for inferring the dynamics of integer- and fractional-Order neuron models. 
Commun. Comput. Phys. 35(1), 2024, 1-37. 
	 
	 \bibitem{Fu}
	 Z.C. Tang, Z.J. Fu,  S. Reutskiy. 
	  An extrinsic approach based on physics-informed neural networks for PDEs on surfaces.  Mathematics. 10(16), 2022, 2861.
	 
	\bibitem{Wang}
 Z.W. Wang, M.X. Chen,  J.R. Chen.
Solving multiscale elliptic problems by sparse radial basis function neural networks.  J. Comput. Phys. 492, 2023, 112452. 
	 
	 \bibitem{Temam}
	  R. T\'{e}mam. Infinite-dimensional dynamical systems in mechanics and physics. New York, Springer-Verlag, 1988.

\bibitem{Tran}
G. Tran, R. Ward. Exact recovery of chaotic systems from highly corrupted data. Multiscale Model. Simul. 15(3),  2007,  1108-1129. 

\bibitem{Turk}
   G. Turk, J. O'brien. Modelling with implicit surfaces that interpolate. ACM Trans. Gragh.  21(4), 2002,  855–873. 
	 
	 \bibitem{Wu19}
	 Z.M. Wu, R. Zhang. Learning physics by data for the motion of a sphere falling in a non-Newtonian fluid. 
	 Commun. Nonlinear Sci. Numer. Simul.  67, 2019, 577-593.
	 
	 \bibitem{Wu_AA}
	 Z.M. Wu, R. Zhang.
	 Probabilistic solutions to DAEs learning from physical data. Anal. Appl. 19(6), 2021, 1093–1111.
	 
	 \bibitem{Zhang}
	 R. Zhang, G. Plonka. Optimal approximation with exponential sums by a maximum likelihood modification of Prony’s method. Adv. Comput. Math. 45, 2019, 1657–1687. 
	 
	 
\end{thebibliography}
\end{document}